\magnification 1000
\input amstex
\documentstyle{amsppt}
\vsize 8.25in
\voffset 1cm
\hoffset 1cm
\topmatter
\rightheadtext{Envelopes of holomorphy and the reflection principle}
\leftheadtext{Jo\"el Merker}
\title
On envelopes of holomorphy of domains with Levi-flat hats
and the reflection principle
\endtitle
\author 
\endauthor
\address 
Laboratoire d'Analyse, Topologie et Probabilit\'es, Centre de
Math\'ematiques et d'Informatique, UMR 6632, 39 rue Joliot Curie,
F-13453 Marseille Cedex 13, France. Fax: 00 33 (0)4 91 11 35 52
\endaddress
\email 
merker\@cmi.univ-mrs.fr 00 33 / (0)4 91 11 36 72 / (0)4 91 53 99 05
\endemail
\thanks
\endthanks

\def\C{{\Bbb C}}
\def\N{{\Bbb N}}
\def\R{{\Bbb R}}
\def\V{{\Bbb V}}

\def\v{\vert}
\def\n{\vert\vert}

\abstract
As expected (or conjectured) after the recent works of
Baouendi-Eben\-felt-Roths\-child, a $\Cal C^\infty$-smooth CR
diffeomorphism $h\: (M,p)\to (M',p')$ between two minimal real
analytic hypersurfaces in $\C^n$ ($n\geq 2$) should be real analytic
{\it if and only if} $(M',p')$ is holomorphically nondegenerate.
Constructing envelopes of holomorphy of special domains, namely ``with
Levi-flat hats'', we provide here a complete original proof of this
statement. As a byproduct of this strategy, we derive an entirely
new treatment of the essentially finite case
(Baouendi-Jacobowitz-Treves' famous theorem, 1985). More generally, we
establish that the reflection mapping $\Cal R_h'$ associated to such a
$\Cal C^\infty$-smooth diffeomorphism between two minimal
hypersurfaces in $\C^n$ ($n\geq 1$) always extends holomorphically to
a neighborhood of $p\times \bar p'$. This gives a generalization of
the Schwarz symmetry principle to higher dimensions.

\medskip










\endabstract

\endtopmatter

\document

\head \S1. Introduction and description of the proof \endhead

The present paper associates the techniques of the reflection
principle and the techniques of analytic discs. Extending CR reflection
objects to a Levi-flat union of Segre varieties, we come down to study the
envelope of holomorphy of certain domains ``with Levi-flat hats''
({\it see} \S2.2 below).

\subhead 1.1.~Main theorem \endsubhead
Let $h\: M\to M'$ be a smooth CR mapping between two real analytic
hypersurfaces in $\C^n$ $(n\geq 2)$, let $p\in M$ and set
$p':=h(p)$. Associated to $h$ and to $M'$ is the so-called {\it
reflection function} (Xiaojun Huang's denomination\,; {\it see}
\S1.6), an interesting invariant, more general than $h$. Our
principal result is\,:

\proclaim{Theorem~1.2}
If $h$ is a $\Cal C^\infty$-smooth CR-diffeomorphism and if $(M,p)$ is
minimal, then the reflection function $\Cal R_h'$ extends
holomorphically to a neighborhood of $p\times \overline{p'}$.  In
particular, $h$ is real analytic at $p$ if $M'$ is holomorphically
nondegenerate.
\endproclaim

\subhead 1.3.~Features of the classical results \endsubhead
The earliest result of this kind was found independently by Hans Lewy
[L] and by Serguei Pinchuk [P1]\,: if $(M,p)$ and $(M',p')$ are
strongly pseudoconvex, then $h$ is real analytic. The classical proof
in [L] and [P1] makes use of the so-called {\it reflection principle}
which consists roughly to solve the map $h$ with respect to $\bar h$
and the jets of $\bar h$ and to apply the Schwarz symmetry principle
in a foliated union of transverse holomorphic discs. Generalizing this
principle, Diederich-Webster proved in 1980 that a sufficiently smooth
CR diffeomorphism is analytic at $p\in M$ if $M$ is generically
Levi-nondegenerate and the morphism of jets of Segre varieties of
$M'$ is injective ({\it see} \S2 in [DW] and \S1.13 below for a
definition).  In 1985, Derridj studied the reflection principle for
proper mappings between some model classes of weakly pseudoconvex
boundaries in $\C^2$. In 1983, Han [Ha] generalized the reflection
principle for CR-diffeomorphisms between what is today called ({\it see}
[BER2]) {\it finitely nondegenerate} hypersurfaces. In 1985, an important
breakthrough was achieved by Baouendi-Jacobowitz-Treves [BJT], who
proved that any $\Cal C^\infty$ CR diffeomorphism $h\: (M,p)\to
(M',p')$ between two real analytic CR-generic manifolds in $\C^n$
which extends holomorphically to a fixed wedge of edge $M$, is real
analytic. After Baouendi-Treves [BT2] (hypersurface case\,; the weakly
pseudoconvex case is treated in [BeFo]), Tumanov [Tu1] (general
codimension), and Baouendi-Rothschild [BR3] (necessity) it was known
that the automatic holomorphic extension to a wedge of the components
of $h$ holds if and only if $(M,p)$ is minimal in the sense of Tumanov
(or equivalently, of finite type in the sense of Bloom-Graham). In the
late eighties, the research on the analyticity of CR mappings has been
pursued by many authors intensively. In 1988, Baouendi-Rothschild and
independently, Diederich-Fornaess extended this kind of reflection principle
to the non diffeomorphic case, namely for a $\Cal C^\infty$ CR map $h$ of
essentially finite
hypersurfaces which is locally finite to one, or locally proper. This
result was generalized in [BR2] to mappings of maximal formal generic
rank on formal Segre varieties (not totally degenerate) or even more
generally with non identically zero formal Jacobian
determinant. Following this circle of ideas, Coupet-Pinchuk-Sukhov 
have pointed out that almost all the above-mentioned reflection
principles come down to the fact that a certain complex analytic
variety $\V_p'$ is zero-dimensional, which intuitively significates
that $h$ is {\it finitely determined by the jets of $\bar h$}, {\it
i.e.} each components $h_j$ of $h$ satisfies a monic Weierstrass
polynomial with some coefficients being analytic functions depending
on a finite jet of $\bar h$ (this observation appears also in
[Me1]). They stated thus a general result in the hypersurface case
whose extension to a higher codimensional minimal CR-generic source
$M$ was achieved recently by Damour in [Da]. In summary, this last
refinement closes up what is attainable in the spirit of the so-called
{\it polynomial identities} devised by Baouendi-Jacobowitz-Treves,
yielding a quite general sufficient condition for the analyticity of
$h$.  In the arbitrary codimensional case, this general sufficient
condition can be expressed as follows. Let $\overline{L}_1,\ldots,
\overline{L}_m$ be a basis of $T^{0,1}M$, denote $\overline{L}^\beta:=
\overline{L}_1^{\beta_1}\ldots \overline{L}_{m}^{\beta_{m}}$ for
$\beta\in\N^{m}$ and let $\rho_{j'}'(t',\bar t')=0$, $1\leq j'\leq d'$, be
a collection of real analytic equations for $(M',p')$. Then the
complex analytic variety
$$
\V_p':=\{t'\in \C^n\: \overline{L}^\beta [\rho'(t',\bar h(\bar p))]=0, \ 
\forall \, \beta \in \N^{m}\}.
\tag 1.4
$$
is always zero-dimensional at $p'\in \V_p'$ in [L], [P1], [DW], [Ha], [De],
[BJT], [BR1], [DF], [BR2], [BR4], [BHR], [DP], [Hu], [BER1], [BER2], 
[CPS1], [CPS2], [Da].  Crucially, the condition $\dim_{p'} \V_p'=0$
requires $(M',p')$ to be essentially finite.

\subhead 1.5.~Non essentially finite hypersurfaces \endsubhead
However, it is known that the finest CR-regularity phenomena come down
to the consideration of a class of much more general hypersurfaces
which are called {\it holomorphically nondegenerate} and which are
generically {\it not essentially finite}. In 1995, 
Baouendi-Rotschild [BR3] exhibited this condition as a {\it necessary and
sufficient} condition for the algebraicity of a biholomorphism between
two real algebraic hypesurfaces. Thanks to the nonlocality of
algebraic objects, they could assume that $(M',p')$ is essentially
finite after a small shift of $p'$ with $\dim_{p'} \V_p'=0$, thus
coming down to known techniques (even in fact simpler, in the
generalization to the higher codimensional case,
Baouendi-Ebenfelt-Rothschild came down to a direct application of the
algebraic implicit function theorem by solving algebraically $h$
with respect to the jets of $\bar h$ [BER1]).  Since then however, few works
have been devoted to the study of the analytic regularity of smooth CR
mapping between non-essentially finite hypersurfaces in $\C^n$. It is
well-known that the main technical difficulties in the subject happen
to occur in $\C^n$ for $n\geq 3$ and that a great deal of such
obstacles can be avoided by assuming that the target hypersurface $M$'
is algebraic (with $M$ real analytic), see {\it e.g.} the works [MM2],
[Mi1,2,3], [CPS1] (in case $M'$ is algebraic, its Segre varieties are
defined all over the compactification $P_{n-1}(\C)$ of $\C^n$, which
helps much). Finally, we would like to mention the papers of Meylan
[Mey], Maire and Meylan [MaMe], Meylan and the author [MM1] in this
concern (nevertheless, after division by a suitable holomorphic
function, the situation under study in these works is again reduced to
polynomial identities).

\subhead 1.6.~Schwarz's reflection principle in higher dimension \endsubhead
In late 1996, seeking a natural generalization of Schwarz's reflection
principle to higher dimension, the author ({\it see} [MM2], [Me1]) discovered
the interest of the so-called {\it reflection function} $\Cal R_h'$
associated with $h$ (the denomination was due to Xiaojun Huang in [Hu], 
but our definition in [Me1] involved one more variable and the crucial 
observation of its biholomorphic invariance) which
appears already implicitely in [BJT]. The explicit expression of this
function depends on a local defining equation for $M'$, but its
holomorphic extendability is independent of coordinates and there are
canonical rules of transformation between two reflection functions.
As the author has pointed out, in the diffeomorphic case, {\it this
function should extend without assuming any nondegeneracy condition on
$M'$, as in the case $n=1$}, provided $M$ is at least of finite type
in the sense of Kohn. It is easy to convince oneself that the
reflection function is the right invariant to study. It has been already
studied thoroughly in the algebraic and in the formal CR-regularity
problems, {\it see} [Me1,3,4], [Mi2,3,4]. For instance, the formal
reflection mapping associated with a formal CR equivalence between two real
analytic CR-generic manifolds in $\C^n$ which are minimal in the sense
of Tumanov is convergent ({\it see} [Mi2,3] for partial
results and [Me4] for the complete statement). If $h$ is a holomorphic
equivalence between two real algebraic CR-generic manifolds in $\C^n$
which are minimal at a Zariski-generic point, then the reflection
mapping $\Cal R_h'$ is algebraic ({\it see} [Mi2] for the
hypersurface case and [Me3] for general codimension).  So we expect that
similar statements hold for smooth mappings between CR manifolds.

\subhead 1.7.~Statement of the results \endsubhead
For our part, we concentrate in this paper on smooth CR mappings
between {\it hypersurfaces}. Thus, let as above $h\: M\to M'$ be a
$\Cal C^\infty$-smooth CR mapping between two small connected pieces
of real analytic hypersurfaces in $\C^n$ with $n\geq 2$. Reserving
generalizations and refinements to further investigation, we shall
assume here for simplicity that $h$ is a CR-diffeomorphism. The
associated reflection function $\Cal R_h'$ is defined over $M\times
\overline{M'}$ as follows. Localizing $M$ and $M'$ at points $p\in M$
and $p'\in M'$ with $p'=h(p)$, we choose a complex analytic defining
equation for $M'$ in the form $\bar z'=\Theta'(\bar w',t')$, where
$t'=(w',z')\in\C^{n-1}\times
\C$ are holomorphic coordinates vanishing at $p'$ and where
the power series $\Theta'(\bar w',t'):=
\sum_{\beta\in\N^{n-1}} (\bar w')^\beta
\, \Theta_\beta'(t')$ vanishes at 
the origin and converges normally in a small polydisc
$\Delta_{2n-1}(0,4r')$, $r'>0$. By definition, the {\it reflection
function} $\Cal R_h'$ associated with $h$ and with such a defining function
is the function $(t,\bar\nu')\mapsto
\bar \mu'-\sum_{\beta\in \N^{n-1}} \bar{\lambda'}^\beta
\, \Theta_\beta'(h(t))=:\Cal R_h'(t,\bar \nu')$, where 
$\bar\nu'=(\bar\lambda',
\bar\mu')\in\C^{n-1}\times \C$. Clearly, this function is CR and smooth
with respect to the variable $t\in M$ and holomorphic with respect to
the variable $\bar\nu'$ near $\overline{M'}$. In case $M$ is minimal
in the sense of Tumanov at a point $p\in M$, then $h$ and $\Cal R_h'$
extend holomorphically to one side $D$ of $M$ at $p$. Our first main
result is as follows.

\proclaim{Theorem~1.8}
If $h$ is a $\Cal C^\infty$-smooth CR-diffeomorphism 
and if $(M,p)$ is minimal, then the reflection function $\Cal R_h'$
extends holomorphically to a neighborhood of $p\times \overline{p'}$.
\endproclaim

\noindent
\remark{Remark}
Of course (exercise), the assumption of minimality of $(M,p)$ can be
swit\-ched to $(M',p')$, because $(M,p)$ and $(M',p')$ are
CR-diffeomorphic.
\endremark

\smallskip
\noindent
Clearly, the holomorphic extendability of $\Cal R_h'$ to a
neighborhood of $p$ is equivalent to the following statement\,: {\it
all the functions $\Theta_\beta'(h(t))=:
\varphi_\beta'(t)$ $($an infinite number$)$ extend holomorphically 
to a neighborhood of $p$} and there exist constants $C_p,r_p>0$ such
that $\v t\v < r_p \Rightarrow \v \varphi_\beta'(t)\v < C_p^{\v
\beta\v+1}$. In certain circumstances, {\it e.g.}  when $(M',p')$ is
Levi-nondegenerate, finitely nondegenerate or essentially finite, one
deduces afterwards that $h$ itself extends holomorphically at $p$. In
Theorem~1.13 below, we shall derive from Theorem~1.8 above an
important expected {\it necessary and sufficient condition} for $h$ to be
holomorphic at $p$.

\subhead 1.9.~Applications \endsubhead We give essentially two 
important applications.

\smallskip
\noindent
$\bullet$ Firstly, associated with $M$', there is an invariant integer
$\kappa_{M'}'$, called the {\it holomorphic degeneracy degree} of
$M'$, which counts the maximal number of holomorphic vector fields
with holomorphic coefficients in a neighborhood of $M'$ which are
tangent to $M'$ and linearly independent at a Zariski-generic point.
We recall that $M'$ is called {\it holomorphically nondegenerate} if
there does not exist a nonzero holomorphic vector field with
holomorphic flow, tangent to $(M',p')$, {\it i.e.} if
$\kappa_{M'}'=0$.  Another (equivalent) definition of $\kappa_{M'}'$
is as follows. Let $j_{t'}^k{\Cal S}_{\tau'}'$ denote the $k$-jet at
the point $t'$ of the complexified Segre variety ${\Cal
S}_{\tau'}'=\{(w',z')\: z'=\bar\Theta'(w',\tau')\}$, which induces a
holomorphic map defined over the extrinsic complexification $\Cal M'$
of $M'$ as follows\,:
$$
\varphi_k' :
{\Cal M}' \ni (t', \tau') \mapsto j_{t'}^k{\Cal
S}_{\tau'}'=(t',\{\partial_{w'}^\beta [z'-\bar\Theta'(w',\tau')]
\}_{\v\beta\v\leq k}) \in \C^{n+{(n-1+k)!\over (n-1)! \, k!}}.
\tag 1.10
$$
We have $\dim_\C \, \Cal M'=2n-1$. Let ${p'}^c:=(p',\bar p')\in\Cal
M'$. It is clear that there exists an integer $\chi_{M'}'$ with $0\leq
\chi_{M'}'\leq n-1$ such that the generic rank of $\varphi_k'$ equals
$n+\chi_{M'}'$ for all $k$ large enough. Then
$\kappa_{M'}'=n-1-\chi_{M'}'$, as is shown in [BER1,2] (for the
algebraic case, {\it see} [Me3]). In general, $M'$ is biholomorphic to
a product $\underline{M}'\times
\Delta^{{\kappa_M'}'}$ by a $\kappa_{M'}'$-dimensional polydisc in a
neighborhood of a Zariski-generic point $q'\in M'$, where
$\underline{M}'\subset \C^{n-\kappa_{M'}'}$ is a {\it holomorphically
nondegenerate} hypersurface. Now, granted Theorem~1.8, we observe that
the graph $\Gamma r (h)=\{(t,h(t))\: t\in (M,p)\}$ of $h$ is clearly
contained in the complex analytic set\,:
$$
\Cal C_h':=\{ (t,t')\in \C^n\times \C^n\:
\Theta_\beta'(t')=\varphi_\beta'(t), \ \forall
\, \beta\in\N^{n-1}\}.
\tag 1.11
$$ 
Since the generic rank of $\varphi_k'$ equals $n+\chi_{M'}'$, there
exists a well-defined irreducible component $\Cal C_h''$ of $\Cal
C_h'$ of dimension $n+\kappa_{M'}'$ containing the graph $\Gamma
r(h)$. We deduce\,:

\proclaim{Corollary~1.12}
Let $\kappa_{M'}'$ be the holomorphic degeneracy degree of $M'$. Then
there exists a pure closed complex analytic subset $\Cal C_h''$ of a
neighborhood of $M\times M'$ in $\C^n$ of dimension $n+\kappa_{M'}'$
which contains the graph of $h$ over $M$. In particular, $h$ extends as 
a correspondence across $M$ if $\kappa_{M'}'=0$.
\endproclaim

\remark{Remark}
According to Coupet-Pinchuk-Sukhov [CPS1,2], a statement equivalent to
Corollary~1.12 above would be that the transcendance degree of the
field extension $\text{\rm Frac}(\Cal O(\Cal V_{\C^n}(M))
\to \text{\rm Frac}(\Cal O(\Cal V_{\C^n}(M)))
(h_1,\ldots,h_n)$ is equal to $\kappa_{M'}'$. It can also be shown
that Corollary~1.12 is equivalent to Theorem~1.8 (exercise).
\endremark

\smallskip
\noindent
$\bullet$ Secondly, an important particular case of Corollary~1.12 is
when $\kappa_{M'}'=0$.  Assuming now that $h$ is $\Cal C^\infty$ in
order to be able to apply a theorem of Malgrange ([Ma],~p.~96]\,: {\it
A $\Cal C^\infty$ manifold is real analytic if and only if it is
contained in a real analytic set of the same dimension}), we deduce
the following important result\,:

\proclaim{Theorem~1.13}
Let $h\: (M,p)\to (M',p')$ be a $\Cal C^\infty$ CR diffeomorphism
between two minimal real analytic hypersurfaces in $\C^n$.  If $M'$ is
holomorphically nondegenerate, then $h$ is real analytic at $p$.
\endproclaim

\noindent
(Of course, real analyticity of $h$ is equivalent to its holomorphic
extendability to a neighborhood of $(M,p)$, by Severi's theorem,
generalized to higher codimension by Tomassini.)

\subhead 1.14.~Necessity \endsubhead
Since 1995-6 ({\it see} [BR3], [BHR]), it is known that Theorems~1.3
above provide an expected {\it necessary and sufficient condition} in
order that $h$ is analytic (provided of course that the local
CR-envelope of holomorphy of $M$, which already contains a side $D$ of
$M$ at $p$, does not contain the other side). Indeed,

\proclaim{Lemma~1.15}
\text{\rm ([BHR])}
Conversely, if $(M',p')$ is holomorphically degenerate and if there
exists a smooth CR function $\varpi\: M'\to \C$ defined in a
neighborhood of $p'\in M'$ \text{\rm which does not extend
holomorphically to a neighborhood of $p'$}, then there exists a
CR-automorphism of $(M',p')$ fixing $p'$ which is not real analytic.
\endproclaim

\subhead 1.16.~Organization\endsubhead
To be brief, \S2 presents first a thorough intuitive description (in
words) of our strategy for the proof of Theorem~1.8, to which the
remainder of the paper is exclusively devoted (since the above applications
are classical).

\subhead 1.17.~Acknowledgement\endsubhead
The author is grateful to Egmont Porten, who pointed out to him the interest
of gluing half-discs to the Levi flat hypersurfaces $\Sigma_\gamma$ below.

\head \S2. Precise description of the proof \endhead

\subhead 2.1.~Envelopes of holomorphy and reflection principle\endsubhead
According to the extendability theorem of Baouendi-Treves [BT2]
(generalized to the $\Cal C^2$-smooth case by Tr\'epreau), the map $h$
in Theorem~1.8 already extends holomorphically to a one-sided
neighborhood $D$ of $M$ at $p$ in $\C^n$. This extension is performed
by using small Bishop discs attached to $M$ and by applying the
Baouendi-Treves approximation theorem [BT1]. By the way, we would like
to remind the reader of the well-known and somewhat paradoxical
phenomenon of {\it automatic holomorphic extension of CR functions on
$M$ to both sides}, which can render the above Theorem~1.8
surprisingly trivial. Indeed, let $U_M$ denote the (open) set of
points $q$ in $M$ such that the envelope of holomorphy of $D$ contains
a neighborhood of $q$ in $\C^n$ (as is well-known, if, for instance,
the Levi form of $M$ has one positive and one negative eigenvalue at
$q$, then $q\in U_M$\,; more generally, the {\it local envelope of
holomorphy} of $M$ or of the one-sided neighborhood 
$D$ of $M$ at an arbitrary point $q\in M$
is always {\it one-sheeted}, as can be proved using the
Baouendi-Treves approximation theorem). Then clearly, the $n$
components of our CR diffeomorphism extend holomorphically to a
neighborhood of $U_M$ in $\C^n$, as does any arbitrary CR function on
$M$. But it remains to extend $h$ holomorphically across $M\backslash
U_M$ and the techniques of the reflection principle are then
unavoidable. {\it Here lies the paradox}\,: sometimes the envelope of
holomorphy trivializes the problem, sometimes it does not
help. Fortunately, in the study of the smooth reflection principle,
the classical techniques do not make usually any difference between
the two sets $U_M$ and $M\backslash U_M$ and these techniques provide
a {\it uniform method} of extending $h$ across $M$, no matter the
reference point $p$ belongs to $U_M$ or to $M\backslash U_M$ ({\it
see} [L], [P1], [DW], [BJT], [BR1], [BR2], [DF], [BHR], [BER1], [BER2],
[CPS1], [CPS2]). Such a uniform method seems to be quite
satisfactory. On the other hand, recent deep works of Pinchuk in the 
study of the {\it geometric reflection principle} show up an accurate 
analysis of the relative pseudo-convex(-concave) loci of $M$. In [P2],
[DP], [Hu], [Sha], the authors achieve the propagation of holomorphic
extension of a germ along the Segre varieties of $M$ (or the Segre
sets), taking into account their relative position with respect to $M$
and its local convexity. In such reasonings, some discussions
concerning envelopes of holomorphy come down naturally in the proofs
(which involve many sub-cases). However, comparing these two trends of
thought, it seems to remain still really paradoxical that both
phenomena contribute to the reflection principle, without an appropriate
understanding of the general links between these two techniques.
Guided by this observation, we have devised a new two-sided technique.
In this article, we shall indeed perform the proof of Theorem~1.8 by
{\it mixing the technique of the reflection principle together with
the consideration of envelopes of holomorphy}. Further, we have been
guided by a deep analogy between the various reflection principles and
the results on {\it propagation of analyticity for CR functions along
CR curves}, in the spirit of the Russian school in the sixties,
of Treves' school, of Tr\'epreau, of Tumanov, of J\"oricke and
others\,: the vector fields of the complex tangent bundle $T^cM$ being
the {\it directions} of propagation for the one-sided holomorphic
extension of CR functions, and the Segre varieties giving these
directions (because $T_q^cM=T_qS_{\bar q}$ for all $q\in M$), one can
expect that they propagate as well the analyticity of CR mappings.  Of
course, such a propagation property is already well-known and
intensively studied since the historical work of Pinchuk
[P1]. However, in the classical works ({\it e.g.} in [P2], [DP]), one
propagates along a single Segre variety $S_{\bar p}$ and perharps
afterwards along the subsequent Segre sets if necessary ({\it see}
[BER1,2], [Me2,4], [Mi3,4], [Sha]). But in this article we will
propagate the analytic properties along a {\it bundle} of Segre
varieties of $M$, namely {\it along a Levi-flat union of Segre
varieties} $\Sigma_\gamma:=\cup_{q\in
\gamma} S_{\bar q}$, parametrized by a smooth curve $\gamma$ transversal to
$T^cM$, in analogy with the propagation of analyticity of CR
functions, where one uses a {\it bundle of attached analytic discs},
parametrized by a curve transversal to $T^cM$.  Here lies the main
displacement of ideas.  Let us now explain our strategy in full
details and describe our proof.

\subhead 2.2.~Description of the proof of Theorem~1.8\endsubhead
To begin with, it is well-knonw that there exists a Zariski-open set
of points $(q',\bar q')$ of $\Cal M'$ at which the rank of the
morphism of jets of Segre varieties 
$\varphi_k'$ is equal to its generic rank, say $n+\chi_{M'}'$,
with $0\leq \chi_{M'}'\leq n-1$. Thus, the jet map $\varphi_k'$ is
locally of constant rank at $(q',\bar q')$. In our first step, we will
show that $\Cal R_h'$ is real analytic at each point $q\in M$ such
that $\varphi_k'$ is locally of constant rank $n+\chi_{M'}'$ at
$(h(q),\overline{h(q)})$. (Of course, using the CR-diffeomorphism
assumption, one can verify that $\chi_M=\chi_{M'}'$, whence also
$\kappa_{M}=\kappa_{M'}'$, but we shall not explicitely need this fact
for the proof.) In fact, if $(M',p')$ is holomorphically nondegenerate,
these points $q'$ are the {\it finitely nondegenerate} points of $M'$,
in the sense of Baouendi-Ebenfelt-Rothschild [BER2]. It will appear
that our proof of the first step, a reminiscence of the Lewy-Pinchuk
reflection principle, appears to be in fact a mild easy generalization
of it. Now, during the second (crucial) step, to which \S4--7 below are
devoted, we shall extend $\Cal R_h'$ across the remaining set
$E_{M'}'$ where $\varphi_k'$ is not locally of constant rank.
This is where we use envelopes of holomorphy. Let
$E_{\text{na}}'\subset E_{M'}'\subset M'$ (``$\text{\rm na}$''
for ``non-analytic'') denote the closed set of points $q'\in M'$ such
that $\Cal R_h'$ is not analytic in a neighborhood of $h^{-1}(q')$. If
$E_{\text{\rm na}}'=\emptyset$, Theorem~1.8 is proved. We shall
therefore assume that $E_{\text{\rm na}}'\neq \emptyset$ and we shall
endeavour to derive a contradiction in several steps as follows.
Following [MP], we shall first show that we can choose a particular
point $p_1'\in E_{\text{\rm na}}'$ which is nicely disposed as follows.

\proclaim{Lemma~2.3}
\text{\rm ({\it cf.} [MP])}
Let $E'\subset M'$ be an arbitrary closed subset of an
everyhere minimal real analytic hypersurface $M'\subset \C^n$, with
$n\geq 2$. If $E'$ and $M'\backslash E'$ are nonempty, 
then there exists a point $p_1'\in E'$ and a $\Cal
C^\omega$ one-codimensional submanifold $M_1'$ of $M'$ with $p_1'\in
M_1'\subset M'$ which is \text{\rm generic} in $\C^n$ and which
divides $M'$ near $p_1'$ in two open parts ${M_1'}^-$ and ${M_1'}^+$
such that $E'$ is contained in the closed side $\overline{{M'}^+}$
near $p_1'$.
\endproclaim

\noindent
To reach the desired contradiction, it will suffice to prove that
$\Cal R_h'$ is analytic at the point $h^{-1}(p_1')$, where $p_1'\in
E_{\text{\rm na}}'\cap M_1'$ is such a point as in Lemma~2.3 above. To
this aim, we shall pick a long embedded real analytic arc $\gamma'$
contained in ${M_1'}^-$ transverse to the complex tangential
directions of $M'$, with the ``center'' $q_1'$ of $\gamma'$ very close
to $p_1'$ and we shall set $E_{\text{\rm na}}:=h^{-1}(E_{\text{\rm na}}')$, 
$\gamma:=h^{-1}(\gamma')$, $p_1:=h^{-1}(p_1')$, $q_1:=h^{-1}(q_1')$,
$M_1:=h^{-1}(M_1')$, $M_1^-=h^{-1}({M_1'}^-)$ and
$M_1^+=h^{-1}({M_1'}^+)$. To the arc $\gamma'$, we shall associate
holomorphic coordinates $t'=(w',z')\in\C^{n-1}\times \C$, $z'=x'+iy'$,
such that $p_1'=0$ and $\gamma'$ is the $x'$-axis (in particular, some
``normal'' coordinates in the sense of Chern-Moser or
Baouendi-Jacobowitz-Treves, called ``regular'' by Ebenfelt, would be
appropriate) and we shall consider the reflection function 
$\Cal R_h'(t,\bar\nu')=
\bar\mu'-\sum_{\beta\in \N^{n-1}} \bar{\lambda'}^\beta \ 
\Theta_\beta'(h(t))$ {\it in these coordinates
$(w',z')$}. The functions $\Theta_\beta'(h(t))$ 
will be called the {\it components of the reflection
function $\Cal R_h'$}. Next, we choose coordinates $t\in \C^n$ near $(M,p_1)$
vanishing at $p_1$. To the $\Cal C^\infty$-smooth arc $\gamma$, we
shall associate the following $\Cal C^\infty$-smooth Levi-flat
hypersurface\,: $\Sigma_\gamma:=
\bigcup_{q\in \gamma} S_{\bar q}$, where $S_{\bar q}$ denotes the 
Segre variety of $M$ associated to various points $q\in M$. Let
$\Delta_n(0,r):=\{t\in \C^n\: \v t\v < r\}$ be the polydisc
with center $0$ of polyradius $(r,\ldots,r)$, where
$r>0$. Using the tangential Cauchy-Riemann operators to
differentiate the fundamental identity which reflects the assumption
$h(M)\subset M'$, we shall establish the following crucial
observation.

\proclaim{Lemma~2.4}
There exists a positive real number $r>0$ independent of
$\gamma'$ such that all the components $\Theta_\beta'(h(t))=\left[
{1\over \beta!}{\partial^{\v
\beta\v}\over \partial {\bar{\lambda'}}^\beta}\Cal R_h'(t,\bar\nu')
\right]\v_{\bar{\lambda'}=0}$ extend as CR functions
of class $\Cal C^\infty$ over $\Sigma_\gamma\cap
\Delta_n(0,r)$.
\endproclaim

\noindent
We now recall that the components $\Theta_\beta'(h(t))$ 
are already holomorphic in $D$ and also
holomorphic in a fixed neighborhood, say $\Omega$, of $M_1^-$ in
$\C^n$, by construction of $M_1'$. In particular, they are holomorphic in
a neighborhood $\omega(\gamma)\subset \Omega$ in $\C^n$ of
$\gamma\subset M_1^-$. Then according to the Hanges-Treves extension
theorem [HaTr], we deduce that all the components
$\Theta_\beta'(h(t))$ of the reflection function extend
holomorphically to a neighborhood $\omega(\Sigma_\gamma)$ of
$\Sigma_\gamma$ in $\C^n$, a (very thin) neighborhood whose size
depends of course on the size of $\omega_\gamma$ (and the size of
$\omega_\gamma$ goes to zero without any explicit control as the
center point $q_1$ of $\gamma$ tends to $p_1\in E_{\text{\rm na}}$).

To achieve the final step, we shall consider the envelope of
holomorphy of $D\cup \Omega \cup \omega(\Sigma_\gamma)$ (in fact, to
prevent from poly-dromy phenomena, we shall instead consider a certain
subdomain of $D\cup \Omega\cup \omega(\Sigma_\gamma)$, see the details
in \S5 below), which is a kind of round domain $D\cup \Omega$ covered
by a thin Levi-flat almost horizontal ``hat-domain''
$\omega(\Sigma_\gamma)$ touching the ``top of the head'' $M$ along the
one-dimensional arc~$\gamma$ (very thin contact). 
Our purpose will be to show that, if the
arc~$\gamma'$ is sufficiently close to $M_1'$ (whence $\gamma$ is also
very close to $M_1$), then the envelope of holomorphy of $D\cup \Omega
\cup \omega(\Sigma_\gamma)$ contains the point $p_1$, {\it
even if $\omega(\Sigma_\gamma)$ is arbitrarily thin}. We will
therefore deduce that all the components of the reflection function
extend holomorphically at $p_1$, thereby deriving the desired
contradiction. By exhibiting a special curved Hartogs domain, we shall
in fact prove that holomorphic functions in $D\cup \Omega \cup
\omega(\Sigma_\gamma)$ extend holomorphically to the lower one sided
neighborhood $\Sigma_\gamma^-$ (roughly speaking, the ``same'' side as
$D=M^-$)\,; we explain below why this analysis gives analyticity at
$p_1$, even in the (simpler) case where $p_1$ belongs to the other
side $\Sigma_\gamma^+$. Notice that, since the order of contact
between $\Sigma_\gamma$ and $M$ is at least equal to two (because
$T_qM=T_q\Sigma_\gamma$ for every point $q\in \gamma$), we cannot
apply directly some version of the edge of the wedge theorem to this
situation. Another possibility (which, on the contrary, works well)
would be to apply repeatedly the Hanges-Treves theorem, in the disc
version given in [Tu2] ({\it see} also [MP]) to deduce that
holomorphic functions in $D\cup \Omega \cup \omega(\Sigma_\gamma)$
extend holomorphically to the lower side $\Sigma_\gamma^-$, just by
sinking progressively $\Sigma_\gamma$ into $D$. But this would require
a too complicated analysis for the desired statement.  Instead, by
performing what seems to be the simplest strategy, we shall use some
deformations (``translations'') of the following half analytic disc
attached to $\Sigma_\gamma$ along $\gamma$. We shall consider the
inverse image by $h$ of the half-disc $(\gamma')^c\cap D'$ obtained by
complexifying $\gamma'$. Rounding off the corners and reparametrizing
the disc, we get an analytic disc $A\in \Cal O(\Delta)\cap
\Cal C^\infty(\overline{\Delta})$ with $A(b^+\Delta)\subset 
\gamma\subset \Sigma_\gamma$, where $b^+\Delta:= b\Delta \cap \{\text{\rm 
Re} \, \zeta \geq 0\}$, $b\Delta=\{\v z\v =1\}$ and $A(1)=q_1$. It is
this half-attached disc that we shall ``translate'' along the complex
tangential directions to~$\Sigma_\gamma$ as follows.

\proclaim{Lemma~2.5}
There exists a $\Cal C^\infty$-smooth 
$(2n-2)$-parameter family of analytic discs $A_\sigma\: 
\Delta\to \C^n$, $\sigma\in\R^{2n-2}$, 
$\v \sigma\v <\varepsilon$, satisfying
\roster
\item"{\bf (1)}"
The disc $A_\sigma\v_{\sigma=0}$ coincides with the above disc $A$.
\item"{\bf (2)}"
The discs $A_\sigma$ are half-attached to $\Sigma_\gamma$, namely
$A_\sigma(b^+\Delta)\subset \Sigma_\gamma$.
\item"{\bf (3)}"
The boundaries $A_\sigma(b\Delta)$ of the discs $A_\sigma$
are contained in $D\cup \Omega \cup \omega(\Sigma_\gamma)$.
\item"{\bf (4)}"
The map $(\zeta,\sigma)\mapsto A_\sigma(\zeta)\in \Sigma_\gamma$ is a 
$\Cal C^\infty$-smooth diffeomorphism from a neighborhood
of $(1,0)\in b\Delta\times\R^{2n-2}$ onto a neighborhood of
$q_1$ in $\Sigma_\gamma$.
\item"{\bf (5)}"
As $\gamma=h^{-1}(\gamma')$ varies and as $q_1$ tends to $p_1$, these 
discs depend $\Cal C^\infty$-smoothly upon $\gamma'$ and properties
{\bf (1-4)} are stable under perturbations of $\gamma'$.
\item"{\bf (6)}"
If $\gamma(0)=q_1$ is sufficiently close to $M_1$, and if $p_1\in
\Sigma_\gamma^-$, then the envelope of
holomorphy of (an appropriate subdomain of) 
$D\cup \Omega \cup \omega(\Sigma_\gamma)$ contains $p_1$.
\endroster
\endproclaim

Consequently, using these properties {\bf (1-6)} and applying the
continuity principle to the family $A_\sigma$, we shall obtain that
the envelope of holomorphy of $D\cup \Omega\cup \omega(\Sigma_\gamma)$
(in fact of a good subdomain of it, in order to assure monodromy)
contains a large part of the side $\Sigma_\gamma^-$ of $\Sigma_\gamma$
in which $D (=:M^-)$ lies. In the case where $p_1$ lies in this side
$\Sigma_\gamma^-$, and provided that the center point $q_1$ of
$\gamma$ is sufficiently close to $p_1$, we are done\,: the components
of the reflection function extend holomorphically at $p_1$. Of course,
it can happen that $p_1$ lies in the side $\Sigma_\gamma^+$ or in
$\Sigma_\gamma$. In fact, the following tri-chotomy is in order to
treat this case. To apply Lemma~2.5 wisely, and to complete the study
of our situation, we shall indeed distinguish three cases.  {\it
Case~I}\,: the Segre variety $S_{\bar p_1}$ cuts $M_1^-$ along an
infinite sequence of points $(q_k)_{k\in\N}$ tending towards $p_1$\,;
{\it Case~II}\,: the Segre variety $S_{\bar p_1}$ does not intersect
$M_1^-$ in a neighborhood of $p_1$ and it goes under $M_1^-$.  {\it
Case~III}\,: the Segre variety $S_{\bar p_1}$ does not intersect
$M_1^-$ in a neighborhood of $p_1$ and it goes over $M_1^-$.  In the
first case, choosing the point $q_1$ above to be one of the points
$q_k$ which is sufficiently close to $p_1$, and using the fact that
$p_1$ {\it belongs to} $S_{\bar q_1}$ (because $q_1\in S_{\bar p_1}$),
we have in this case $p_1\in \Sigma_\gamma$ and the holomorphic
extension to a neighborhood $\omega(\Sigma_\gamma)$ already yields
analyticity at $p_1$ (in this case, we have nevertheless to use
Lemma~2.5 to insure monodromy of the extension). In the second case,
we have $S_{\bar p_1}\cap D\neq\emptyset$. We then choose the center
point $q_1$ of $\gamma$ very close to $p_1$. Because we then have a
uniform control of the size of $\omega(\Sigma_\gamma)$, we again get
that $p_1$ belongs to $\omega(\Sigma_\gamma)$ and Lemma~2.5 is again
used to insure monodromy. In the third (a priori more delicate) case, by a
simple calculation, we shall observe that $p_1$ always belong to
$\Sigma_\gamma^-$ and Lemma~2.5 applies to yield holomorphic extension
and monodromy of the extension, and we are done in the three cases. In
conclusion to this presentation, we would like to say that some
unavoidable technicalities that we have not mentioned here will render the
proof a little bit more complicated (especially about the choice of
$q_1$ sufficiently close to $p_1$, about the choice of $\gamma$ and
about the smooth dependence with respect to $\gamma$ of
$\Sigma_\gamma$ and of $A_\sigma$).

\head \S3. Extension across a Zariski dense open subset of $M$ \endhead

The starting point of the proof is to show that the 
the reflection function extends holomorphically
across $M$ at a Zariski-generic point of $M$.
This is done by performing a very easy generalization of the
classical Lewy-Pinchuk reflection principle. Thus, let
$\varphi_k'$ denote the morphism of $k$-jets of Segre varieties of
$M'$ expressed in some coordinate system as in~\thetag{1.10} above. The
generic rank of the holomorphic map $\varphi_k'$ stabilizes for $k$
large enough. We let $n+\chi_{M'}'$, $0\leq \chi_{M'}'\leq n-1$,
denote this generic rank and we fix such a large integer $k$
definitely. We set $\kappa':=\kappa_{M'}'=n-1-\chi_{M'}'$. We denote
by $E_{M'}'$ the set of points $p'\in M'$ in a neighborhood of which
the Segre morphism $\varphi_k'$ is not of constant rank $n+\chi_{M'}'$.
Using the representation~\thetag{1.10} of $\varphi_k'$ in coordinates,
it can be checked there exists a complex analytic set $\Cal E'$,
closed in a neighborhood of $M'$, such that $E_{M'}'=\Cal E'\cap M'$
(exercise\,; {\it cf.}~[Me1]\,; we believe that $\Cal E'$ is the
intrinsic complexification of $E_{M'}'$, but we have no proof of it\,;
anyway, we need not such a result). We set $G_{M'}':= M'\backslash
E_{M'}'$. The following fact is well known in the subject ({\it see} {\it
e.g.} [BER2], [Me3]).

\proclaim{Lemma~3.1}
For each point $q'\in G_{M'}'$, there exists a local coordinate system
$(w',v',z')\in \C^{n-1-\kappa'}\times
\C^{\kappa'}\times \C$ vanishing at $q'$ such that the defining equations
of $(M',q')$ take the simple form $\bar z'=\Theta'(\bar w',w',z')$,
where the analytic function $\Theta'$ is independent
of the coordinates $(v',\bar v')$.
\endproclaim

\noindent
In other words, near such a Zariski-generic point $q'\in G_{M'}'$,
then $M'$ is biholomorphic to the product $\underline{M}'\times
\Delta^{\kappa'}$, where the hypersurface $\underline{M}'\subset
\C^{n-\kappa'}$, equipped with coordinates $(w',z')$, is given by the
equation $\bar z'=\Theta'(\bar w',w',z')$. As $q'\in G_{M'}'$, we have
furthermore that the morphism $\underline{\varphi}_k'$ of $k$-jets of
Segre varieties of $\underline{M}'$ is immersive at $\underline{q}'$,
{\it i.e.}  has rank maximal equal to $2(n-\kappa')-1$. In other
words, the hypersurface $\underline{M}'\subset \C^{n-\kappa'}$ is
finitely nondegenerate at $\underline{q}'$ ([BER2]).  The main result
of this paragraph can be summarized as follows. It then shows that
$\Cal R_h'$ is analytic at each point of $G_{M'}'$, since $M$ is
minimal at $q$ for all $q$ in a small neighborhood of $p$. In the
following assertion, we use $\Cal C^1$-smoothness only (as noticed in
[BJT], [HMM], in the $\Cal C^1$ case, one can prove more generally
that the components of the reflection function extend holomorphically
across $p$ in the essentially finite case\,; we prove the elementary
Lemma~3.2 just for completeness).

\proclaim{Lemma~3.2}
Let $h\: (M,q)\to (M',q')$ be a $\Cal C^1$-smooth CR-diffeomorphism. If
$(M,q)$ is minimal and the Segre morphism of $M'$ is of constant maximal rank
at $q'$ $(${\it i.e.} $q'\in G_{M'})$, then $\Cal R_h'$ extends 
holomorphically at $q$.
\endproclaim

\demo{Proof}
We choose coordinates $(w',v',z')$ as above, we split the components
of the mapping as $h=(g,l,f)$ accordingly and we set
$\underline{h}:=(g,f)$. The reflection function $\Cal
R_h'=\bar\mu'-\sum_{\beta\in \N^{n-1-\kappa_{M'}'}}
\bar{\lambda'}^\beta \ \Theta_\beta'(\underline{h})$
is independent of the $\kappa'$ components
$l=(l_1,\ldots,l_{\kappa'})$ of the mapping $h$ and $\Cal R_h'$
depends only on the partial map $\underline{h}$. We shall therefore
immediately deduce that $\Cal R_h'$ is analytic at $0$, from the
following assertion. \qed\enddemo

\proclaim{Lemma~3.3}
The components of $\underline{h}=(g,f)$ extend holomorphically at $q$.
\endproclaim

\demo{Proof}
The proof is an easy generalization of the Lewy-Pinchuk reflection
principle, which corresponds to the case $\kappa'=0$ and $M'$ being
$2$-nondegenerate at $q'$ in the sense of [BER2].  Let
$L_1,\ldots,L_{n-1}$ be a commuting basis of $T^{1,0}M$ with real
analytic coefficients, for instance $L_j={\partial \over \partial
w_j}+\bar\Theta_{w_j}(w,\bar t){\partial
\over \partial z}$, $j=1,\ldots,n-1$, where $z=\bar\Theta(w,\bar t)$ is a
defining equation for $(M,q)$ with coordinates $t=(w,z)$ vanishing at
$q$. Since $h$ is a $\Cal C^1$-smooth CR diffeomorphism, after an
eventual linear change of coordinates near $M$, we can assume that the
determinant $\hbox{det} \, \left(\bar L_j \bar g_k(\bar
t)\right)_{1\leq j,k\leq n-1-\kappa'}$ is nonzero at the point $\bar
t=0$. Applying the derivations $\overline{L}_1,\ldots, 
\overline{L}_{n-1-\kappa'}$ to the
fundamental identity $\bar f=\Theta'(\bar g,\underline{h})$ and using
Cramer's rule as in [P1],
we deduce that for each multiindex $\beta\in
\N^{n-1-\kappa'}$ with $\v\beta\v=1$, then there exists an analytic
function $\underline{\Omega}_\beta$ such that, for all $t\in M$ in a
neighborhood of $0$, we have\,:
$$
[\partial_{\bar w'}^\beta\Theta'](
\bar g(t),\underline{h}(t))=\underline{\Omega}_\beta(t,
\bar t,\{\partial_t^\alpha
\overline{\underline{h}}(t)\}_{\v\alpha\v \leq 1})=:
\underline{\omega}_\beta(t,\bar t).
\tag 3.4
$$
Here, the right-hand sides of~\thetag{3.4} are {\it a priori} only
$\Cal C^0$-smooth with respect to $t\in M$, but we observe readily 
that the left-hand
sides are in fact $\Cal C^1$. Therefore the right-hand sides
$\underline{\omega}_\beta(t,\bar t)$ are in fact $\Cal C^1$ over
$M$. Thus, we shall be allowed to apply again the de\-rivations
$\overline{L}_1,\ldots,\overline{L}_{n-1-\kappa'}$ to the 
equations~\thetag{3.4}. To begin with,
let $D=M^-$ denote the local one-sided neighborhood of $M$ to which the
components of $h$ have a holomorphic extension, by the
Baouendi-treves extension theorem [BT2] (we have assumed that
$(M,q)$ is minimal) and let $M^+$ denote the other side, approximatively 
symmetric to $D$. As in the Lewy-Pinchuk reflection
principle, using the one-dimensional Schwarz reflection in the complex
lines $\{w=ct.\}$ which are transverse to $M$, we observe that the
functions $\underline{\omega}_\beta$ extend continuously to $M^+$ as
functions which are partially holomorphic with respect to the
transverse variable $z$. Since their boundary value 
$\underline{\omega}_\beta$ is $\Cal C^1$ on the boundary $(M,q)$
thanks to the above observation, a known regularity principle in
one-dimensional complex analysis ({\it see} [H\"o], [BJT], [HMM])
shows that the partial holomorphic
extension with respect to $z$ of $\underline{\omega}_\beta$ to $M^+$
is in fact of class $\Cal C^1$ over $M^+\cup M$. 
We can thus re-apply the derivations $\overline{L}_j$'s to
the identity $ [\partial_{\bar w'}^\beta\Theta'](
\bar g(t),\underline{h}(t))=\underline{\omega}_\beta(t,\bar t)$ and then
use again Cramer's rule, to deduce that for each $\beta\in \N^{n-1-\kappa'}$
with $\v\beta\v =2$, there exists a $\Cal C^0$ function
$\underline{\omega}_\beta(t,\bar t)$, extending holomorphically with
respect to $z$ in $M^+$ such that $[\partial_{\bar w'}^\beta\Theta'](
\bar g(t),\underline{h}(t))=\underline{\omega}_\beta(t,\bar t)$. Again, 
we deduce from this relation that such $\underline{\omega}_\beta$'s
for $\v \beta \v=2$ are in fact $\Cal C^1$ over $M^+\cup M$ and we get in
conclusion a similar identity by induction on $\beta$ for all
$\beta\in \N^{n-1-\kappa'}$. Let us write this relation in the form\,:
$$
\Theta_\beta'(\underline{h}(t))+\sum
{\,}_{\gamma\in\N_*^{\nu}} \
\bar g(\bar t)^\gamma \ \Theta_{\beta+\gamma}'(\underline{h}(t)) \
(\beta+\gamma)!/[\beta! \, \gamma!]=
 \ \underline{\omega}_\beta(t,\bar t)/[\beta!],
\tag 3.5
$$
where $\nu:=n-1-\kappa'$ and $\N_*^\nu:=\N^\nu\backslash \{0\}$.  Now,
it is easy to check that the map $\underline{t}'\mapsto
(\Theta_\beta'(\underline{t}'))_{\v\beta\v \leq k}$ is immersive at
$0$, since the map $(\underline{t}',\underline{\tau}')\mapsto
(\underline{\tau}',([\partial_{\bar w'}^\beta \Theta'](\bar
w',\underline{t}'))_{\v
\beta\v \leq k})$ is
immersive at $0$ on $\underline{\Cal M'}$, by assumption, for $k$ 
large enough. In
eqs.~\thetag{3.5}, we consider the terms $\bar g(\bar t)^\gamma$ with
$\bar g(0)=0$ to be $\Cal C^1$ over $M^+\cup M$ and partially
holomorphic with respect to $z$, as are the left-hand side terms
$\underline{\omega}_\beta(t,\bar t)$.  Applying therefore the implicit
function theorem to eqs.~\thetag{3.5}, we deduce that there exists a
$\Cal C^1$ mapping $a=(a_1,\ldots, a_{n-\kappa'})$ over $M$ which
extend partially holomorphically with respect to $z$ into $M^+$ such
that $\underline{h}(t)=a(t,\bar t)$ when $t\in M$. As in the classical
Lewy-Pinchuk reflection principle, this proves that $\underline{h}$ extends
holomorphically at $0$. The proofs of Lemmas~3.2 and~3.3 are complete
now.
\qed
\enddemo

\subhead 3.6.~On equivalences of hypersurfaces\endsubhead
As an application of Lemma~3.2, we mention here the following easy and
very useful corollary (see \S7 below).

\proclaim{Lemma~3.7}
Let $h\: (M,p)\to (M',p')$ be a $\Cal C^k$ CR-diffeomorphism, 
$1\leq k<\infty$. If
$(M,p)$ is minimal, if the Segre morphism of $M'$ is of constant rank
in a neighborhood of $(p',\bar p')$ and if $M'$ is given by an arbitrary
equation of the form $\bar z'=\sum_{\beta\in\N^{n-1}} \bar{w'}^\beta \
\Theta_\beta'(t')$, then there exists a holomorphic map $H\: (M,p)\to
(M',p')$ whose $k$-th jet at $p$ coincides with the $k$-jet at $p$ of
$h$ such that $\Theta_\beta'(H(t))\equiv
\Theta_\beta'(h(t))$ for all $\beta\in\N^{n-1}$.
Furthermore, if $\kappa_{M'}'=0$, then $H$ is unique, $H\equiv h$ and
$h$ is analytic. Finally, for $q$ running in a neighborhood of $p$,
there exists such a family of equivalences $H_q\,: (M,q)\to (M',q')$
depending $\Cal C^k$-smoothly with respect to $q$.
\endproclaim

\remark{Remark}
We believe that Lemma~3.7 holds true without the restriction that
the Segre morphism is locally of constant rank, but we have no proof
of that (on the other hand, if $h$ is $\Cal C^\infty$, we have a proof
of it, just by considering the formal map induced by the Taylor series
of $h$ at $p$ and by applying the ``Corollaire~2.7'' of [Me4]).
\endremark

\demo{Proof}
The biholomorphic invariance of Segre varieties entails that this
property $\Theta_\beta'(H(t))\equiv
\Theta_\beta'(h(t))$ for all $\beta\in\N^{n-1}$ 
is satisfied for every system of coordinates if and only if it is
satisfied for a single such system of coordinates. Now, in the
coordinates given by Lemma~3.1, we have established that $h=(g,l,f)$
has the property that the components $g$ and $f$ are analytic. It
suffices therefore to choose $H_q:=(g,j_q^kl,f)$ for $q$ in a
neighborhood of $p$.
\qed
\enddemo

\head \S4. Layout of a typical point of non analyticity \endhead

Thus, we already know that $\Cal R_h'$ is analytic over the open dense
subset $h^{-1}(G_{M'}')$ of $M$. We notice that using the
considerations of \S2 above, it can be shown that
$G_M=h^{-1}(G_{M'}')$, but we shall not need this fact. It remains to
show that $\Cal R_h'$ is analytic at each point $p:=h^{-1}(p')$ for
$p'\in E_{M'}'$. This objective constitutes the principal task of the
demonstration. In fact, we shall prove a slightly more general
semi-global statement which we can summarize as follows.

\proclaim{Theorem~4.1}
Let $h\: (M,p)\to (M',p')$ be a $\Cal C^\infty$-smooth
CR-diffeomorphism of connected everywhere
minimal real analytic hypersurfaces in $\C^n$. If the local reflection 
mapping $\Cal R_h'$ associated to 
some coordinate system for $(M',p')$
is analytic at one point $q$ of $(M,p)$, then it is 
analytic all over $(M,p)$.
\endproclaim

\remark{Remark}
If the small piece of hypersurface $(M,p)$ is minimal at $p$, then
shrinking it if necessary, it will be minimal at every point, because
the condition $\hbox{Lie}_q (T^cM)=T_qM$ is an open condition.
Thus, of course we assume that $(M,p)$ is 
minimal at every point and in fact, 
our Theorem~4.1 holds true if, instead of a germ,
we consider a map of {\it globally minimal} (hence connected) ``large''
real analytic hypersurfaces.
\endremark

\smallskip

\subhead 4.2.~Construction of a generic wall \endsubhead
In \S2 above, we have already shown that $\Cal R_h'$ is analytic at
each point of the nonempty open set $h^{-1}(G_{M'}')$, thus
Theorem~4.1 implies our goal, Theorem~1.8. To prove Theorem~4.1, let us
denote by $E_{\text{\rm na}}'$ the closed set of points $p'\in M'$
such that $\Cal R_h'$ is not analytic in a neighborhood of
$h^{-1}(p')$. If $E_{\text{\rm na}}'=\emptyset$, we are done. We
suppose therefore by contradiction that $E_{\text{\rm na}}'\neq
\emptyset$. We shall reach a contradiction by showing that there
exists {\it a} ({\it i.e.} {\it at least one}) {\it point} 
$p_1'\in E_{\text{\rm na}}'$
such that $\Cal R_h'$ is analytic at $h^{-1}(p_1')$. As in Lemma~2.3,
this point $p_1'$ will belong to a generic one-codimensional
submanifold $M_1'\subset M'$, a kind of ``wall'' in $M'$ dividing $M'$
locally into two open sides, which will be disposed conveniently in
order that one open side of the ``wall'', say ${M_1'}^-$, will {\it
contain only points where $\Cal R_h'$ is already real analytic}. To
show the existence of such a point $p_1'\in E_{\text{\rm na}}'$ and
of such a manifold (``wall'') $M_1'$, we shall proceed as in [MP,~Lemma~2.3] 
(the reader is referred to this paper
for more details about this construction).
As $M'$ is minimal at every point, it is certainly globally
minimal, in the sense that the CR-orbit of every point of $M'$ is $M'$
itself in the whole. We thus choose an arbitrary point $r'\in
M'\backslash E_{\text{\rm na}}'\neq \emptyset$ (it is nonempty,
because we already know that 
$M'\backslash E_{\text{\rm na}}'\supset M'\backslash
E_{M'}'=G_{M'}'\neq \emptyset$). The CR-orbit of $r'$ is equal to $M'$
in the whole. Thus there exists a piecewise differentiable $\Cal
C^\omega$ curve $\gamma'$ running in complex tangential directions to
$M'$ with origin $r'$ and endpoint a point $q'\in E_{\text{\rm na}}'$.
After shortening $\gamma'$, we can assume that $\gamma'$ is a smoothly
embedded real segment with boundary with $q'=E_{\text{\rm na}}'\cap
\gamma'$ and that $\gamma'$ extends a little bit further as the
integral curve of a $\Cal C^\omega$ section $L'$ of $T^cM'$. Then
using the dynamical flow of $L'$ and a one parameter family of real
$(2n-2)$-dimensional spheroids $B_s$, $s\in \R$, which are stretched
along the flow lines of $L'$ in $M'$ and centered at a point $r'\in
\gamma'\backslash E_{\text{\rm na}}'$ sufficiently close to $q'$ ({\it
see} [MP,~Lemma~2.3]), we can choose $M_1'$ to be a piece of the
boundary of the first spheroid $B_{s_1}$ which happens to touch
$E_{\text{\rm na}}'$ at a point $p_1'$. Furthermore, the
$(2n-2)$-dimensional spheroid $B_{s_1}$ is {\it generic} in $\C^n$ at
$p_1'$ by construction, since we can assume in this construction that
the vector field $L'$ is tangent to the spheroids $B_s$ {\it only}
along a fixed $(2n-3)$-dimensional spheroid contained in $M'\backslash
E_{\text{\rm na}}'$ which is independent of $s$. In summary, it
suffices now for our purposes to establish the following assertion.

\proclaim{Theorem~4.3}
Let $p_1'\in E_{\text{\rm na}}'$ and assume that there exists a $\Cal
C^\omega$ one-codimensional submanifold $M_1'$ with $p_1'\in
M_1'\subset M'$ which is generic in $\C^n$ such that $E_{\text{\rm na}}'
\backslash \{p_1'\}$ is completely contained in one of the 
two open sides of
$M'$ divided by $M_1'$ at $p_1'$, say in ${M_1'}^+$, and such that 
$\Cal R_h'$ is analytic at each point 
of the other side ${M_1'}^-$. Then the reflection function $\Cal R_h'$
extends holomorphically at $h^{-1}(p_1')$.
\endproclaim

\noindent
Now, an elementary reasoning using only linear changes of coordinates
and Taylor's formula shows that, after an eventual change of the
manifold $M_1'$ in a new manifold $M_1''$ which is bent quadratically
in the side ${M_1'}^-$, we can assume that $M'$ is given by the
equation $z'=\bar\Theta'(w',\bar t')$, that $M_1'$ is given by the two
equations $z'=\bar\Theta'(w',\bar t')$, $u_1'=-[{v_1'}^2+\v w_*'\bar
w_*'\v^2+x^2]$, where $w_1'=u_1'+iv_1'$, $w_*'=(w_2',\ldots,
w_{n-1}')$ and that the side ${M_1'}^-$ is given by\,:
$$
{M_1'}^-\: \ \ \ \{(w',z')\in M'\: \ u_1'<-[{v_1'}^2+\v w_*'\v^2+x^2]\},
\tag 4.4
$$
in coordinates $(w',z')$ centered at $p_1'$. We set
$p_1:=h^{-1}(p_1')$, $M_1:=h^{-1}(M_1')$ which is a $\Cal
C^\infty$-smooth one-codimensional generic submanifold of $M$. By
assumption, the reflection function $\Cal R_h'$ associated with these
coordinates is already holomorphic at each point of the side of
automatic extension $D=\{(w,z)\in \C^n\: y<h(w,\bar w,x)\}$. {\it and
it is also real analytic at each point of $M_1^{-}\subset M$}. Let us
write this more precisely. Without loss of generality, we can assume
that $\Theta$ and $\Theta'$ both converge in $\Delta_n(0,4r)$, with
$r>0$. Let $(\Psi_{p'}')_{p'\in M'}$ denote a family of
biholomorphisms sending $p'\in M'$ to $0$, holomorphically
parametrized by $p'\in \Delta_n(0,r)$ with $\Psi_0'=\text{\rm Id}$.
Let $\bar z'=\Theta_{p'}'(\bar w',t')$ denote the equation of
$\Psi_{p'}'(M')$. By saying that $\Cal R_h'$ (associated with $p'=0$)
extends holomorphically across $M$ at each point of $M_1^-$, we mean
precisely that each reflection function $\Cal R_{p',h_{p'}}$ in these
coordinates extends holomorphically to a neighborhood of $p\times 0$,
for every point $p\in M_1^-$. Using then an explicit representation
for $\Psi_{p'}'$ and achieving elementary calculations with power
series, we obtain the following concrete characterization in which a
single coordinate system is considered.

\proclaim{Lemma~4.5}
There exists a neighborhood $\Omega$ of $M_1^-$ and a lower semi-continuous
positive function $r(p)>0$ such that for each point $p\in M_1^-$, the
polydisc $\Delta_n(p,r(p))$ is contained in $\Omega$, the components
$\Theta_\beta'(h(t))$ of the reflection function associated with a
fixed coordinate system converge in
$\Delta_n(0,r(p))$ and they satisfy a Cauchy estimate $\v
\Theta_\beta'(h(t))\v < C_p^{\v \beta\v +1}$ for $\v t-p\v < r(p)$.
\endproclaim

\noindent
Now start our principal constructions. As explained in \S2.2, we
intend to study the envelope of holomorphy of the union of $D$
together with an arbitrary thin neighborhood of a Levi-flat
hypersurface $\Sigma_\gamma$. We need real arcs and analytic discs.

\head \S5. Envelopes of holomorphy of domains with Levi-flat hats \endhead

\subhead 5.1.~A family of real analytic arcs \endsubhead
To begin with, we choose coordinates $t$ and $t'$ as above near $M$
and near $M'$ in which $p_1$ and $p_1'$ are the origin and in which
the equations of $M$ and of $M'$ are given by
$M\,: z=\bar\Theta(w,\bar t)$ and $M'\: z'=\bar\Theta'(w',\bar t')$. We
can assume that the power series defining $\bar\Theta$ and $\bar
\Theta'$ converge normally in the polydisc $\Delta_{2n-1}(0,4r)$, for
some $r>0$. In the target space, we 
now define a convenient, sufficiently rich, family of
embedded real analytic arcs $\gamma_{w_{q'}',x_{q'}'}'(s')$, depending on
$(2n-1)$ very small real parameters $(w_{q'}',x_{q'}')\in \C^{n-1}
\times \R$ satisfying $\v w_{q'}'\v < \varepsilon$, 
$\v x_{q'}'\v < \varepsilon$, where $\varepsilon << r$, with ``time'' $s'$
satisfying $\v s'\v \leq 2r$,
and which are all transverse to the complex tangential 
directions of $M'$ (shrinking $r$ if necessary), as follows\,:
$$
\left\{
\aligned
&
\gamma_{w_{q'}',x_{q'}'}':=\left\{(u_{1,q'}'-{s'}^2-(v_{1,q'}'+s')^2-
\v w_{*q'}'\v^2-{x_{q'}'}^2+i[v_{1,q'}'+s'],\right.\\
&
\ \ \ \ \ \ \ \ \ \ \ \ \ \ \ \ \ \ \ \ \ \ \ \ \ \ \ \ \ \
\ \ \ \ \ \ \ \ \ \ \ \ \ \ \
\left.
,w_{*q'}', \, x_{q'}') \in M'\,:
\ s'\in \R, \, \v s' \v \leq 2r\right\}.
\endaligned\right.
\tag 5.2
$$
It can be straightforwardly checked that the 
following properties hold\,:
\roster
\item"{\bf (1)}"
{\it The mapping $(w_{q'}',x_{q'}')\mapsto \gamma_{w_{q'}',x_{q'}'}'(0)$
is a $\Cal C^\omega$ real diffeomorphism onto a neighborhood of 
$0$ in $M'$. Furthermore, the inverse image of $M_1'$ and of ${M_1'}^-$
correspond to the sets $\{u_{1,q'}'=0\}$ and $\{u_{1,q'}'<0\}$,
respectively.}
\item"{\bf (2)}"
{\it For $u_{1,q'}'<0$, we have $\gamma_{w_{q'}',x_{q'}'}'\subset \subset
{M_1'}^-$.}
\item"{\bf (3)}"
{\it For $u_{1,q'}'=0$, we have 
$\gamma_{w_{q'}',x_{q'}'}'\cap M_1'=\{\gamma_{w_{q'}',x_{q'}'}'(0)\}$.}
\item"{\bf (4)}"
{\it The order of contact of $\gamma_{w_{q'}',x_{q'}'}'$ with 
$M_1'$ at the point $\gamma_{w_{q'}',x_{q'}'}'(0)$ equals $2$.}
\endroster 

\subhead 5.3.~Inverse images\endsubhead 
Since $h$ is a $\Cal C^\infty$
CR diffeomorphism, we get in $M$ a family of $\Cal
C^\infty$ arcs, $h^{-1}(\gamma_{w_{q'}',x_{q'}'}')$. It is clear
that we get a parameterized family of arcs depending in a $\Cal
C^\infty$ fashion with respect to some variables $(w_q,x_q,s)\in
\C^{n-1}\times \R\times \R$, which we will denote by
$\gamma_{w_q,x_q}(s)$ accordingly (by the index notation
$\cdot_{w_q,x_q}$, we mean that the arc is parametrized by its center
point $\gamma_{w_q,x_q}(0)\in M$, which covers in a diffeomorphic way a
neighborhood of $0$ in the manifold $M$ equipped with coordinates
$(w,x)$\,: this is why we shall maintain in the sequel such a
notation). Of course, after adapting a bit the domains of variation
(or shrinking a bit $\varepsilon$ and $r$), we can suppose that $\v
w_q\v < \varepsilon$, $\v x_q\v < \varepsilon$, $\v s\v
\leq 2r$ and again $\varepsilon << r$. Then the $\Cal C^\infty$ arcs
$\gamma_{w_q,x_q}$ satisfy the two properties {\bf (1)}, {\bf (2)},
{\bf (3)} and {\bf (4)}
above with respect to $M_1$. In particular,
\roster
\item"{\bf (5)}"
{\it There exists a continuous function $r(\varepsilon)$ with
$0<r(\varepsilon)\leq 2 r$ and tending 
to $0$ with $\varepsilon$ such that, for all
$(w_q,x_q)$ with $\v w_q\v, \v x_q\v < \varepsilon$, we have}\,:
$$
\{\gamma_{w_q,x_q}(s)\: \ r(\varepsilon)\leq s \leq 2 r\}
\subset \subset M_1^-.
\tag 5.4
$$
\endroster
This property will be of interest later, when envelopes appear on scene. 

\subhead 5.5.~Construction of a family of Levi-flat hats\endsubhead
Next, if $\gamma$ is a $\Cal C^\infty$-smooth arc in $M$ transverse to
$T^cM$ at each point, we can construct the union of Segre varieties
attached to the points running in $\gamma$\,: $\Sigma_{\gamma}:=
\bigcup_{p\in \gamma} S_{\bar p}$. For various arcs $\gamma_{w_q,x_q}$,
we obtain various sets $\Sigma_{\gamma_{w_q,x_q}}$ which are in fact $\Cal
C^\infty$-smooth Levi-flat hypersurfaces in a neighborhood of
$\gamma_{w_q,x_q}$. The uniformity of the size of such neighborhoods
follows immediately from the smooth dependence with respect to
$(w_q,x_q)$. What we shall need in the sequel can be then summarized
as follows.

\proclaim{Lemma~5.6}
After shrinking $r$ if necessary, there exists $\varepsilon>0$
with $\varepsilon << r$ such that, if the parameters
$(w_q,x_q)$ satisfy $\v w_q\v, \v x_q\v < \varepsilon$, then
the set $\Sigma_{\gamma_{w_q,x_q}}$ is a closed $\Cal C^\infty$-smooth
$($and $\Cal C^\infty$-smoothly parametrized$)$ 
\text{\rm Levi-flat} hypersurface of
$\Delta_n(0,r)$.
\endproclaim

\subhead 5.7.~Two families of half-attached analytic discs\endsubhead
Let us now define inverse images of analytic discs. Complexifying
the $\Cal C^\omega$ arcs $\gamma_{w_{q'}',x_{q'}'}'$, we obtain
transverse holomorphic discs, closed in $\Delta_n(0,3r/2)$, of which
one half part penetrates inside $D':=h(D)$. Uniformly smoothing out
the corners of such half discs, using Riemann's conformal mapping
theorem and then an automorphism of $\Delta$, we can easily construct
a $\Cal C^\omega$-parameterized family of analytic discs
$A_{w_{q'}',x_{q'}'}'\:
\Delta\to \C^n$ which are $\Cal C^\infty$ up to the boundary $b\Delta$
such that, if we denote $b^+\Delta:= b\Delta\cap \{\text{\rm Re} \,
\zeta \geq 0\}$ (and $b^-\Delta:= b\Delta\cap \{\text{\rm Re} \,
\zeta \leq 0\}$), then we have $A_{w_{q'}',x_{q'}'}'(1)=
\gamma_{w_{q'}',x_{q'}'}'(0)$ and also\,:
$$
A_{w_{q'}',x_{q'}'}'(b^+\Delta) \subset \gamma_{w_{q'}',x_{q'}'}'
\ \ \ \ \ \text{\rm and} \ \ \ \ \
A_{w_{q'}',x_{q'}'}'(b^+\Delta) \supset \gamma_{w_{q'}',x_{q'}'}'
\cap \Delta_n(0,5r/4),
\tag 5.8
$$
for all $\v w_{q'}'\v, \v x_{q'}'\v < \varepsilon$. Consequently, 
the composition with 
$h^{-1}$ yields a family of analytic discs $A_{w_q,x_q}:=
h^{-1}\circ A_{w_{q'}',x_{q'}'}'$ which are half-attached to 
$\Sigma_{\gamma_{w_q,x_q}}$ and which 
satisfy similar properties, namely\,:
\roster
\item"{\bf (1)}"
{\it The map $(w_q,x_q,\zeta)\mapsto A_{w_q,x_q}(\zeta)$ is
$\Cal C^\infty$-smooth and it provides a uniform
family of $\Cal C^\infty$ embeddings 
of $\overline{\Delta}$ into $\C^n$.}
\item"{\bf (2)}"
{\it We have $A_{w_q,x_q}(1)=\gamma_{w_q,x_q}(0)$, and also}\,:
$$
A_{w_q,x_q}(b^+\Delta) \subset \gamma_{w_q,x_q}
\ \ \ \ \ \text{\rm and} \ \ \ \ \
A_{w_q,x_q}(b^+\Delta) \supset \gamma_{w_q,x_q}\cap \Delta_n(0,r).
\tag 5.9
$$
\item"{\bf (3)}"
{\it $A_{w_q,x_q}(b^-\Delta) \subset\subset D\cup M_1^-$.}
\endroster
This family $A_{w_q,x_q}$ will be our starting point to study the
envelope of holomorphy of (a certain subdomain of) the union of 
$D$ together with a neighborhood $\Omega$ of
$M_1^-$ and an arbitrarily thin neighborhood of $\Sigma_{w_q,x_q}$.
At first, we must include $A_{w_q,x_q}$ 
into a larger family of discs obtained
by sliding the half-attached part inside $\Sigma_{\gamma_{w_q,x_q}}$
along its complex tangential directions.

\subhead 5.10.~Deformation of half-attached analytic discs\endsubhead
To this aim, we introduce the equation $y=H_{w_q,x_q}(w,x)$
of $\Sigma_{\gamma_{w_q,x_q}}$, where the map $(w_q,x_q,w,x)\mapsto
H_{w_q,x_q}(w,x)$ is of course $\Cal C^\infty$ and 
$\n H_{w_q,x_q}-H_{0,0}\n_{C^\infty(w,x)}$ is very small.
Further, we need some formal notation. We denote $A_{w_q,x_q}(\zeta):=
(w_{w_q,x_q}(\zeta),z_{w_q,x_q}(\zeta))$ and
$A_{w_q,x_q}(1)[=\gamma_{w_q,x_q}(0)]:=
(w_{w_q,x_q}^1,z_{w_q,x_q}^1)$. For our discs, we shall choose the
regularity $\Cal C^{1,\alpha}$, $0<\alpha<1$, which is sufficient for
our purposes.  Let $T_1$ denote the Hilbert transform vanishing at $1$
({\it see} [Tu], [MP], [BER2]\,; by definition, $T_1$ is the unique
(bounded, by a classical result) endomorphism of $\Cal
C^{1,\alpha}(\overline{\Delta},\R)$, $0<\alpha<1$, to itself such that
$\phi+iT_1(\phi)$ extends holomorphically to $\Delta$ and $T_1\phi$
vanishes at $1\in b\Delta$, {\it i.e.}  $(T_1(\phi))(1)=0$). Let
$\varphi^-$ and $\varphi^+$ be two $\Cal C^\infty$ functions on
$b\Delta$ satisfying $\varphi^-\equiv 0$, $\varphi^+\equiv 1$ on
$b^+\Delta$ and $\varphi^-+\varphi^+=1$ on 
$b\Delta$. The fact that our discs are half attached to
$\Sigma_{\gamma_{w_q,x_q}}$ can be expressed by saying that
$y_{w_q,x_q}(\zeta)=\varphi^+(\zeta) \, H_{w_q,x_q}(
w_{w_q,x_q}(\zeta),x_{w_q,x_q}(\zeta))+\varphi^-(\zeta) \
y_{w_q,x_q}(\zeta)$ for all $\zeta\in b\Delta$ 
({\it cf.} A\u\i rapetyan [A\u\i]). Since the two
functions $x_{w_q,x_q}$ and $y_{w_q,x_q}$ on $b\Delta$ are harmonic
conjugates, the following (Bishop) equation is satisfied on
$b\Delta$ by $x_{w_q,x_q}$\,:
$$
x_{w_q,x_q}(\zeta)=-[T_1 (\varphi^+ \, H_{w_q,x_q}(w_{w_q,x_q},\, 
x_{w_q,x_q}))](\zeta) +\psi_{w_q,x_q}(\zeta)+
x_{w_q,x_q}^1, 
\tag 5.11
$$
where we have set $\psi_{w_q,x_q}(\zeta):=-[T_1(\varphi^-
\, y_{w_q,x_q})](\zeta)$. We want to perturb these
discs $A_{w_q,x_q}$ by ``translating'' them along the complex
tangential directions to $\Sigma_{\gamma_{w_q,x_q}}$. Introducing a
new parameter $\sigma\in\C^{n-1}$ with $\v \sigma \v < \varepsilon$,
we can now include the discs $A_{w_q,x_q}$ into a larger
parameterized family $A_{w_q,x_q,\sigma}$ by solving the following 
perturbed Bishop equation on $b\Delta$ with parameters $(w_q,x_q,\sigma)$\,:
$$
x_{w_q,x_q,\sigma}(\zeta)=-
[T_1 (\varphi^+ \,
H_{w_q,x_q}(w_{w_q,x_q}+\sigma\, ,x_{w_q,x_q,\sigma}))](\zeta)+
\psi_{w_q,x_q}(\zeta)+x_{w_q,x_q}^1.
\tag 5.12
$$
For instance, the existence and the $\Cal C^{1,\beta}$-smoothness
(with $0<\beta<\alpha$ arbitrary) of a solution $x_{w_q,x_q,\sigma}$
to~\thetag{5.12} follows from Tumanov's work [Tu3]. Clearly the
solution disc $A_{w_q,x_q,\sigma}$ is half attached to 
$\Sigma_{\gamma_{w_q,x_q}}$. Differentiating the Bishop 
equation~\thetag{5.12} with respect to $\sigma$, one sees that 
the derivatives $(\partial / \partial \sigma)(z_{w_q,x_q,\sigma})$
is uniformly small ({\it cf.}~a similar computation in [Tu1,2,3]). 
In summary\,:

\proclaim{Lemma 5.13}
After shrinking perharps $\varepsilon$, there exists a $\Cal
C^{1,\beta}$-smooth mapping defined for $\v w_q\v, \v x_q\v <
\varepsilon$ and for $\sigma\in \C^{n-1}$, 
$\v \sigma \v < \varepsilon$, $(w_q,x_q,\sigma,\zeta)\mapsto
A_{w_q,x_q,\sigma}(\zeta)$, which is holomorphic with respect to
$\zeta$, and which fulfills the following conditions\,:
\roster
\item"{\bf (1)}"
$A_{w_q,x_q,0}\equiv A_{w_q,x_q}$.
\item"{\bf (2)}"
$A_{w_q,x_q,\sigma}(b^+\Delta)\subset \Sigma_{\gamma_{w_q,x_q}}$ for all
$\sigma$.
\item"{\bf (3)}"
The map $\C^{n-1}\times b^+\Delta\ni (\sigma,\zeta) 
\mapsto A_{0,0,\sigma}(\zeta)\in \Sigma_{\gamma_{0,0}}$ 
is a local $\Cal C^{1,\beta}$
diffeomorphism from a neighborhood of $0\times 1$ onto
a neighborhood of $A_{0,0}(1)=0$.
\endroster
\endproclaim

\subhead 5.14.~Preliminary to the continuity principle \endsubhead
We are now in position to state and to prove the main assertion of this
paragraph. At first, we shall let the parameters $(w_q,x_q,\sigma)$
range in certain new precise subdomains. We choose a positive
$\delta < \varepsilon$ with the property that the range of the map in
{\bf (3)} above, when restricted to $\{\v\sigma\v < \delta\}\times
b^+\Delta$, {\it contains} the intersection of $\Sigma_{\gamma_{0,0}}$
with a small polydisc $\Delta_n(0,2\eta)$, for some 
$\eta>0$. Of course, there exists a
constant $c>1$, depending only on the Jacobian of the map in {\bf (3)}
at $0\times 1$ such that ${1\over c}
\delta \leq \eta \leq c \delta$. Furthermore, since
the boundary of the disc $A_{0,0,0}$ is transversal to $T_0^c
\Sigma_{\gamma_{0,0}}$ (whence ${d\over d\lambda}
A_{0,0,0}(\lambda)\v_{\lambda\in \R,
\lambda=1}\not\in T_0M)$, then after shrinking a bit
$\eta$ if necessary, we can assume that the set
$\{A_{0,0,\sigma}(\zeta)\,:
\v \sigma \v < \delta, \zeta\in \Delta\cap \Delta(1,\delta)\}$ {\it contains
and foliates by half analytic discs the whole lower side} 
$\Delta_n(0,2\eta)\cap \Sigma_{\gamma_{0,0}}^-$.

\remark{Remark}
Of course, the side $\Sigma_{\gamma_{0,0}}^-$ is ``the same side''
as $M^-$, {\it i.e.} the side of $\Sigma_{\gamma_{0,0}}$ 
where the greatest portion of $D$ lies, although $D$ is in
general {\it not} entirely contained in $\Sigma_{\gamma_{0,0}}^-$,
because the Segre varieties $S_{\bar q}$ for $q\in \gamma_{0,0}$
may well intersect $D$, as is known.
\endremark

\smallskip
\noindent
As in \S2, we now fix a neighborhood $\Omega$ of $M_1^-$ in $\C^n$ to
which all the components of the reflection function extend
holomorphically as in Lemma~4.5. We have already shown that the half parts
$A_{w_q,x_q,\sigma}(b^+\Delta)$ are all contained in
$\Sigma_{\gamma_{w_q,x_q}}$ (hence in arbitrarily thin neighborhoods
of it). It remains now to control the half parts
$A_{w_q,x_q,\sigma}(b^-\Delta)$. Using property {\bf (3)}
after~\thetag{5.9} for $(w_q,x_q)=(0,0)$, namely
$A_{0,0,0}(b^-\Delta)\subset\subset D\cup M_1^-$, it is clear that,
after shrinking $\delta$ if necessary, then we can insure that
$A_{0,0,\sigma}(b^-\Delta) \subset \subset D\cup \Omega$ for all
$\sigma\in \C^{n-1}$ with $\v \sigma \v < \delta$. Of course, this
shrinking will result in a simultaneous shrinking of $\eta$, and we
still have the important supclusion\,:
$\{A_{0,0,\sigma}(\zeta)\,:
\v \sigma \v < \delta, \zeta\in \Delta\cap \Delta(1,\delta)\} \supset
\Delta_n(0,2\eta)\cap \Sigma_{\gamma_{0,0}}^-$. Finally, shrinking again
$\varepsilon$ if necessary, we then come to a situation that we may 
summarize\,:
\roster
\item"{\bf (1)}"
{\it For all $\v w_q\v, \v x_q\v < \varepsilon$, we have}\,:
$$
\{A_{w_q,x_q,\sigma}(\zeta)\,:
\v \sigma \v < \delta, \zeta\in \Delta\cap \Delta(1,\delta)\}\supset
\Delta_n(0,\eta)\cap \Sigma_{\gamma_{w_q,x_q}}^-.
\tag 5.15
$$
\item"{\bf (2)}"
{\it $A_{w_q,x_q,\sigma}(b^+\Delta)\subset \Sigma_{\gamma_{w_q,x_q}}$ and 
$A_{w_q,x_q,\sigma}(b^-\Delta)\subset \subset D\cup \Omega$, for all
$\v \sigma \v < \delta$}.
\endroster

\remark{Remark}
Shrinking for the last time $\delta$ and $\varepsilon$ if necessary,
we can further insure that all the dics $A_{w_q,x_q,\sigma}$ are
embeddings of $\overline{\Delta}$ in $\C^n$, which will be convenient to
apply the continuity principle. If $\varepsilon$
is small enough, we can also insure that the
intersection of $D$ with $\Delta_n(0,\eta)\cap
\Sigma_{\gamma_{w_q,x_q}}^-$ is {\it connected} for all $\v w_q\v, \v x_q\v
< \varepsilon$.
\endremark

\subhead 5.16.~Envelopes of holomorphy \endsubhead
We are now in position to state and to prove the main assertion of this
paragraph. Especially, the following lemma will be applied to 
each member of the collection $\{\Theta_\beta'(h(t))\}_{\beta\in\N^{n-1}}$.

\proclaim{Lemma~5.17}
Let $\delta,\, \eta,\, \varepsilon>0$ be as above. If a holomorphic
function $\psi\in\Cal O(D)$ extends holomorphically to a neighborhood
$\omega(\Sigma_{\gamma_{w_q,x_q}})$, then there exists a unique
holomorphic function $\Psi\in \Cal O(D\cup [\Delta_n(0,\eta)\cap
\Sigma_{\gamma_{w_q,x_q}}^-])$ such that $\Psi\v_D\equiv \psi$.
\endproclaim

\demo{Proof}
Here, we fix $w_q$ and $x_q$. First, we notice that the assumption
entails that $\psi$ extends to a holomorphic function in domains of
the form $D_1\cup \omega(\Sigma_{\gamma_{w_q,x_q}})$, for any
subdomain $D_1\subset D$ is such that $D_1\cap
\omega(\Sigma_{\gamma_{w_q,x_q}})$ is {\it connected}. Clearly, $D_1$
can be chosen so that $D_1\cup \omega(\Sigma_{\gamma_{w_q,x_q}})$
contains the union of the boundaries $A_{w_q,x_q,\sigma}(b\Delta)$,
for $\v \sigma \v < \delta$. Next, we notice that for various
$\sigma$'s, all our discs are clearly {\it analytically isotopic} to
$A_{w_q,x_q,0}$ ({\it see} [Me2]). Since this last disc is clearly
isotopic to a point in $D_1$ (just do the isotopy by shrinking
$A_{w_q,x_q}$ to its center point $A_{w_q,x_q}(0)$), we can apply the
continuity principle in the version given by Lemma~3.2 in [Me2] to
deduce that, for all $\sigma$, there exists a holomorphic function in
a neighborhood of $A_{w_q,x_q,\sigma}(\overline{\Delta})$ which
coincides with $\psi$ in a neighborhood of
$A_{w_q,x_q,\sigma}(b\Delta)$. Then some arguments similar to those
given in [Me2,~p.~40] apply to deduce that there exists a unique
holomorphic function $\chi \in \Cal O(\Delta_n(0,\eta)\cap
\Sigma_{\gamma_{w_q,x_q}}^-)$ which coincides with $\psi$ in a 
neighborhood of $\gamma_{w_q,x_q}$ (the fact that the map
$\Delta_{n-1}(0,\delta)\times \Delta \cap \Delta(1,\delta) \ni
(\sigma,\zeta)\mapsto A_{w_q,x_q,\sigma}(\zeta)\in \C^n$ is an
embedding is an important property which insures uniqueness of the extension).
Finally, using the connectedness of $\Delta_n(0,\eta)\cap
\Sigma_{\gamma_{w_q,x_q}}^-\cap D$ 
and the principle of analytic continuation,
we get the holomorphic function 
$\Psi\in\Cal O(D\cup [\Delta_n(0,\eta)\cap
\Sigma_{\gamma_{w_q,x_q}}^-])$, 
which completes the proof of Lemma~5.17.
\qed
\enddemo

\head \S6. Holomorphic extension to a Levi-flat union of Segre varieties
\endhead

\subhead 6.1.~Straightenings \endsubhead
For each parameter $(w_{q'}',x_{q'}')$, we have considered
the analytic arc $\gamma_{w_{q'}',x_{q'}'}'$ defined by~\thetag{5.2}.
To this family of analytic arcs we can clearly associate a 
family of normalizing coordinates as follows.

\proclaim{Lemma~6.2} 
If $\varepsilon<<r$ is small enough, there exists a $\Cal
C^\omega$-parameterized family of biholomorphic mappings
$\Phi_{w_{q'}',x_{q'}'}'$ of $\Delta_n(0,2r)$ which straightens
$\gamma_{w_{q'}',x_{q'}'}'$ to the $x'$-axis, such that the image
$M_{w_{q'}',x_{q'}'}':=\Phi_{w_{q'}',x_{q'}'}(M')$ is a closed
$\C^\omega$ hypersurface of $\Delta_n(0,r)$ close to $M_{0,0}'$ in
$\C^\omega$ norm which is given by an equation of the form
$\bar z'=\Theta_{w_{q'}',x_{q'}'}'(\bar w',t')$, with
$\Theta_{w_{q'}',x_{q'}'}'(\bar w',t')$ converging normally in the
polydisc $\Delta_{2n-1}(0,r)$ and satisfying
$\Theta_{w_{q'}',x_{q'}'}'(0,t')\equiv 0$.
\endproclaim

\subhead 6.3.~Different reflection functions \endsubhead
Let us develope this equation in the form\,:
$$
\bar z'=\Theta_{w_{q'}',x_{q'}'}'(\bar w',t')=\sum_{\beta\in \N^{n-1}}
\bar{w'}^\beta \ \Theta_{w_{q'}',x_{q'}',\beta}'(t').
\tag 6.4
$$
We denote by
$h_{w_{q'}',x_{q'}'}=(g_{w_{q'}',x_{q'}'},f_{w_{q'}',x_{q'}'})$ the
map in these coordinate systems. To all such coordinates are
therefore associated {\it different reflection functions} by\,:
$$
\Cal R_{w_{q'}',x_{q'}',h_{w_{q'}',x_{q'}'}}'(t,\bar\nu')
:=\bar\mu'-\sum_{\beta\in\N^{n-1}}
\bar{\lambda'}^\beta \ \Theta_{w_{q'}',x_{q'}',\beta}'
(h_{w_{q'}',x_{q'}'}(t)).
\tag 6.5
$$ 

\subhead 6.6.~Holomorphic extension to a Levi-flat hat \endsubhead
We now come to the main construction of this paragraph. Let
$E_{M'}'=\Cal E'\cap M'\subset M'$ be the real analytic subset of \S3.
Then we claim that $E_{M'}'$ is of real dimension $\leq 2n-3$. Indeed,
the complex dimension of $\Cal E'$ is $\leq n-1$. If $E_{M'}'$ would
contain a $(2n-2)$-dimensional real analytic piece of manifold, this
piece would certainly be a complex analytic hypersurface contained in
$M'$, contradicting the fact that $M'$ is minimal at every
point. Because the codimension of $E_{M'}'$ in $M'$ is $\geq 2$, then for
almost all $(w_{q'}',x_{q'}')$ ({\it i.e.} except a closed set of zero
Lebesgue measure in the parameter space), then the intersection
$\gamma_{w_{q'}',x_{q'}'}'\cap E_{M'}'$ is empty. In this 
situation, we have\,: 

\proclaim{Lemma~6.7}
For every $(w_{q'}',x_{q'}')$ with $\gamma_{w_{q'}',x_{q'}'}'\cap
E_{M'}'=\emptyset$ and $\gamma_{w_{q'}',x_{q'}'}'\subset
{M_1'}^-$, then all the components
$\Theta_{w_{q'}',x_{q'}',\beta}'(h_{w_{q'}',x_{q'}'}(t))$ of $\Cal R_{w_{q'}',x_{q'}',h_{w_{q'}',x_{q'}'}}'
(t,\bar\nu')$ extend to be holomorphic in a neighborhood
$\omega(\Sigma_{\gamma_{w_q,x_q}})$ of $\Sigma_{\gamma_{w_q,x_q}}$ in
$\C^n$.
\endproclaim

\remark{Remark}
To prove Lemma~6.7, we need the existence of a map $H$ as in Lemma~3.7
{\it at every point of $\gamma_{w_q,x_q}$}. This is why we require
$\gamma_{w_q,x_q}\cap E_{M'}'\neq\emptyset$. However, as $h$ is $\Cal
C^\infty$ in Theorem~1.8, the condition $\gamma_{w_{q'}',x_{q'}'}'\cap
E_{M'}'=\emptyset$ can be removed, thanks to the remark after
Lemma~3.7. But we have in mind the same theorem with $h$ being only of
class $\Cal C^{1,\alpha}$ ({\it see} Theorem~7.12 below).  Thus, we
will conduct the proof taking into account $E_{M'}'$, even when $h$ is
$\Cal C^\infty$.
\endremark

\demo{Proof}
After a biholomorphic change of coordinates near $M$, we can assume
from the beginning that $dh(0)=\hbox{Id}$ and that $\Theta(\bar w,t)$
converges normally in $\Delta_{2n-1}(0,4r)$. At first, we shall
establish our main crucial observation as follows.

\proclaim{Lemma~6.8}
If $\v w_q\v, \v x_q\v < \varepsilon$ and $\varepsilon$ is
sufficiently small, then all the components $\Theta_\beta'(h(t))$
extend as CR functions of class $\Cal C^\infty$ over
$\Sigma_{\gamma_{w_q,x_q}}\cap
\Delta_n(0,r)$.
\endproclaim

\noindent
Now we remind the reader that the $\Theta_\beta'(h(t))$'s already
extend holomorphically to a neighborhood
$\omega(\gamma_{w_q,x_q})\subset\Omega$ of $\gamma_{w_q,x_q}\subset
M_1^-$ in $\C^n$, by construction. Taking Lemma~6.8 for granted we shall
then complete the proof of Lemma~6.7 by an application of the
following statement (a minor generalization of the Hanges-Treves
extension theorem to parametrized hypersurfaces which also holds
in $\Cal C^{1,\alpha}$)\,: \qed \enddemo

\proclaim{Lemma~6.9}
Let $\Sigma$ be a $\Cal C^\infty$-smooth Levi-flat hypersurface in $\C^n$
$(n\geq 2)$. If a continuous CR function $\psi$ extends
holomorphically to a neighborhood $\Cal U_p$ of a point $p$ belonging
to a leaf $\Cal F_\Sigma$ of $\Sigma$, then $\psi$ extends
holomorphically to a neighborhood $\omega(\Cal F_\Sigma)$ of $\Cal
F_\Sigma$ in $\C^n$. The size of this neighborhood $\omega(\Cal
F_\Sigma)$ depends on the size of $\Cal U_p$ and is stable under
sufficiently small $($even non-Levi-flat$)$ perturbations of $\Sigma$.
\endproclaim

\demo{Proof of Lemma~6.8} 
Let $\overline{L}_1,\ldots,\overline{L}_{n-1}$
be the commuting basis of $T^{0,1}M$ given by $\overline{L}_j={\partial
\over \partial \bar w_j} +\Theta_{\bar
w_j}(\bar w,t){\partial \over \partial \bar z}$, $1\leq j\leq
n-1$. Clearly, the coefficients of these vectors fields converge
normally in the polydisc $\Delta_{2n-1}(0,4r)$. By the diffeomorphism
assumption, we have $\hbox{det} (\overline{L}_j \bar g_k)_{1\leq j,k\leq
n-1}(0)\neq 0$. We shall denote this determinant by\,:
$$
\underline{\Cal D}(\bar w,t,\{\partial_{\bar t_j}
\bar g_k(\bar t) \}_{1\leq j\leq n, \, 1\leq k\leq n-1}).
\tag 6.10
$$ 
Here, $t$ belongs to $M$ and the function $\underline{\Cal D}$ is
holomorphic in its variables. Replacing $z$ by $\bar\Theta(w,\bar t)$
in $\underline{\Cal D}$, we can write $\underline{\Cal D}$ in the form
$\Cal D(w,\bar t,\{\partial_{\bar t_j} \bar g_k(\bar t) \}_{1\leq j\leq n,
\, 1\leq k\leq n-1})$, where $\Cal D$ is holomorphic in its variables. 
Shrinking $r>0$ if necessary, we may assume that for all fixed
coordinate point $\bar t_q=(\bar w_q,\bar z_q)\in M$ with $\v \bar
t_q\v<r$, then\,:
\roster
\item"{\bf (1)}"
{\it The polarization $\Cal D(w,\bar t_q,\{\partial_{\bar t_j}
\bar g_k(\bar t_q) \}_{1\leq j\leq n, \, 1\leq k\leq n-1})$ is convergent
on the Segre variety $S_{\bar t_q}\cap \Delta_n(0,2r)=\{(w,z)\in
\Delta_n(0,2r)\: z=\bar\Theta(w,\bar t_q)\}$, {\it i.e.} is convergent
with respect to $w$ for $\v w \v < 2r$}.
\item"{\bf (2)}" 
{\it This expression $\Cal D(w,\bar t_q,\{\partial_{\bar t_j}
\bar g_k(\bar t_q) \}_{1\leq j\leq n, \, 1\leq k\leq n-1})$ 
does not vanish at any point of $S_{\bar t_q}\cap
\Delta_n(0,2r)$, {\it i.e.} does not vanish for all $\v w \v < 2r$}.
\endroster 
Let us choose $(w_{q'}',x_{q'}')$ satisfying
$\gamma_{w_{q'}',x_{q'}'}'\subset {M_1'}^-$, with $
\v w_{q'}'\v, \v x_{q'}'\v <\varepsilon$. We pick the corresponding 
parameter $(w_q,x_q)$ with $\v w_q\v, \v x_q\v <
\varepsilon$. By the choice of $\Phi_{w_{q'}',x_{q'}'}'$, we then have
$g_{w_{q'}',x_{q'}'}(\gamma_{w_q,x_q}(s))=0$ for all $s\in \R$ with
$\v s\v \leq 2r$. This property will be crucial. As $h$ is only of
class $\Cal C^\infty$ over $M$, to apply the tangential Cauchy-Riemann
derivations and to make a polarization afterwards, we need at first
replace $h_{w_{q'}',x_{q'}'}$ in a neighborhood of
$\gamma_{w_q,x_q}(s)$ by a local {\it holomorphic} equivalence
$H_{w_{q'}',x_{q'}'}\: (M,\gamma_{w_q,x_q}(s))\to (M',
\gamma_{w_{q'}',x_{q'}'}'(s'))$ given by Lemma~3.7 (with $k=1$) satisfying\,:
\roster
\item"{\bf (1)}"
$H_{w_{q'}',x_{q'}',s}(\gamma_{w_q,x_q}(s))=
h_{w_{q'}',x_{q'}'}(\gamma_{w_q,x_q}(s))=
(0,f_{w_{q'}',x_{q'}'}(\gamma_{w_q,x_q}(s)))$.
\item"{\bf (2)}"
$dH_{w_{q'}',x_{q'}'}(\gamma_{w_q,x_q}(s))=
dh_{w_{q'}',x_{q'}'}(\gamma_{w_q,x_q}(s))$.
\item"{\bf (3)}" 
$\Theta_{w_{q'}',x_{q'}',\beta}'
(H_{w_{q'}',x_{q'}'}(t))\equiv 
\Theta_{w_{q'}',x_{q'}',\beta}'
(h_{w_{q'}',x_{q'}'}(t))$ for all $t$ 
close to $\gamma_{w_q,x_q}(s)$.
\endroster

\noindent
For this application of Lemma~3.7, we fix the parameter $s$, namely we
fix the point $\gamma_{w_q,x_q}(s)$ (but according to Lemma~3.7, when
$s$ varies, for all $k\in\N$, there will exist such holomorphic
equivalence $H_{w_{q'}',x_{q'}',s}$ depending $\Cal C^k$-smoothly with
respect to $s$, a property which we shall need below).  If we write
$H_{w_{q'}',x_{q'}'}:=(G_{w_{q'}',x_{q'}'},F_{w_{q'}', x_{q'}'})$,
then applying the tangential Cauchy-Riemann operators
$\overline{L}_1^{\beta_1}\cdots
\overline{L}_{n-1}^{\beta_{n-1}}$, $\beta\in \N^{n-1}$ to the identity\,:
$$
\bar F_{w_{q'}',x_{q'}'}(\bar t)=\Theta_{w_{q'}',x_{q'}'}'(\bar
G_{w_{q'}',x_{q'}'}(\bar t), H_{w_{q'}',x_{q'}'}(t)),
\tag 6.11
$$ 
which holds for $t$ in a neighborhood of $\gamma_{w_q,x_q}(s)$, we get
by a classical calculation ({\it see} [BJT], [BR1,4], 
[BER2], [Me4]) an infinite family
of identities of the following kind, for all $\beta\in\N^{n-1}$
(the case $\beta=0$ simply means~\thetag{6.11})\,:
$$
\left\{
\aligned
&
\Theta_{w_{q'}',x_{q'}',\beta}'(H(t))+\sum_{\kappa\in\N_*^{n-1}}
{(\beta+\kappa)!\over \beta!\ \kappa!} \ \bar G(\bar t)^\gamma
\ \Theta_{w_{q'}',x_{q'}',\beta+\kappa}'(H(t))=
\\
& 
\ \ \ \ \ \ \ \ \ \ \ \ \ \ \ \ \ \ \ \ \
={\Cal T_\beta(w,\bar t,\{\partial_{\bar t}^\kappa 
\bar H_j(\bar t)\}_{1\leq j\leq n, \, \v
\kappa\v \leq \v \beta \v})\over 
[\Cal D(w,\bar t,\{\partial_{\bar t_j}
\bar G_k(\bar t) \}_{1\leq j\leq n, 
\, 1\leq k\leq n-1})]^{2\v \beta\v -1}}.
\endaligned\right.
\tag 6.12
$$
Here the $\Cal T_\beta$'s are holomorphic with respect to $(w,\bar t)$
and relatively polynomial with respect to the jets
$\{\partial_{\bar t_j} \bar G_k(\bar t) \}_{1\leq j\leq n,
\, 1\leq k\leq n-1}$. Also, the variable $t$ runs in $M$ in a neighborhood 
of $\gamma_{w_q,x_q}(s)$. To lighten a bit the notation, we have
dropped the subscript of $H_{w_{q'}',x_{q'}'}$ exceptionally (however,
we will anyway keep this subscript in the sequel).  Now, we remind the
reader that all the functions
$t\mapsto\Theta_{w_{q'}',x_{q'}',\beta}'(h_{w_{q'}',x_{q'}',
\beta}'(t))$ are already real analytic in a neighborhood of 
$\gamma_{w_q,x_q}$, since $\gamma_{w_q,x_q}\subset M_1^-$. Let us
denote them by\,:
$$
\varphi_{w_{q'}',x_{q'}',\beta}'(t):= \Theta_{w_{q'}',x_{q'}',\beta}'
(h_{w_{q'}',x_{q'}'}(t))\equiv 
\Theta_{w_{q'}',x_{q'}',\beta}'
(H_{w_{q'}',x_{q'}'}(t)).
\tag 6.13
$$ 
We can replace them directly in the left hand side 
of~\thetag{6.12}. Written in this form, ~\thetag{6.12} then involves
only functions which are holomorphic in $t$ and in $\bar t$, for
$(t,\bar t)$ running in a neighborhood of $(\gamma_{w_q,x_q}(s),
\overline{\gamma_{w_q,x_q}(s)})$. Thus, we can complexify~\thetag{6.12},
by replacing $(t,\bar t)$ by $(t,\tau)\in \Cal M$ close to
$(\gamma_{w_q,x_q}(s), \overline{\gamma_{w_q,x_q}(s)})$, where $\Cal
M$ is the extrinsic complexification of $M$, given in coordinates
$(t,\tau)=(w,z,\zeta,\xi)$ by the holomorphic equation $\xi
=\Theta(\zeta,t)$. Choosing $(t,\tau)$ of the form
$(w,z,\overline{\gamma_{w_q,x_q}(s)})$, namely, choosing $(w,z)$ to
belong to the Segre variety $S_{\overline{\gamma_{w_q,x_q}(s)}}$, and
using the important fact that $\bar
g(\overline{\gamma_{w_q,x_q}(s)})=\bar
G(\overline{\gamma_{w_q,x_q}(s)})=0$ for such a fixed $s$ (this is a
crucial point, since it entails that the queue sum
$\sum_{\kappa\in\N_*^{n-1}}$ in~\thetag{6.12} disappears), we obtain,
for $t\in S_{\overline{\gamma_{w_q,x_q}(s)}}$ close to
$\gamma_{w_q,x_q}(s)$, {\it i.e.} for 
$t=(w,\bar\Theta(w,\overline{\gamma_{w_q,x_q}(s)}))$\,:
$$
\varphi_{w_{q'}',x_{q'}',\beta}'(w,\bar\Theta(w,
\overline{\gamma_{w_q,x_q}(s)}))
={\Cal T_\beta(w,\bar t,
\{\partial_{\bar t}^\alpha \bar H_k(\bar t)\}_{1\leq
k\leq n, \ \v \alpha\v \leq
\v \beta \v})\over 
[\Cal D(w,\bar t,\{\partial_{\bar t_j}
\bar g_k(\bar t) \}_{1\leq j\leq n, 
\, 1\leq k\leq n-1})]^{2\v \beta\v -1}}.
\tag 6.14
$$
Here, we have directly replaced $\partial_{\bar t_j}
\bar G_k(\bar t)$ by $\partial_{\bar t_j} \bar g_k(\bar t)$ 
in the denominator, which is allowed, since the one-jets of $\bar
H_{w_{q'}',x_{q'}'}$ and of $\bar h_{w_{q'}',x_{q'}'}$ coincide at
$\overline{\gamma_{w_q,x_q}(s)}$. Thanks to {\bf (1)} and {\bf (2)}
after~\thetag{6.10}, we see that the right hand side of~\thetag{6.14}
converges with respect to $w$ for all $\v w \v < r$. We thus have got
that the Taylor series of the components of the reflection function
converge on each leaf of $\Sigma_{\gamma_{w_q,x_q}}\cap
\Delta_n(0,r)$. It remains finally to verify that the right hand sides
of~\thetag{6.14} depend in fact in a $\Cal C^\infty$ fashion with
respect to $s$. Applying Lemma~3.7, we see that we can construct
convergent equivalences $H_{w_{q'}',x_{q'}',s}$ depending $\Cal
C^k$-smoothly with respect to $s$. Replacing it in~\thetag{6.14}, we
get the desired transversal smooth dependence over
$\Sigma_{\gamma_{w_q,x_q}}$.  This completes the proofs of Lemmas~6.7
and~6.8.
\qed
\enddemo

\remark{Remark}
As $h$ is $\Cal C^\infty$ in our Theorems~1.8 and~1.13, it is in
fact superfluous to replace $h_{w_{q'}',x_{q'}'}$ by a holomorphic
local equivalence $H_{w_{q'}',x_{q'}',s}$ at $\gamma_{w_q,x_q}(s)$,
because infinite derivations are allowed and we see in this case
that~\thetag{6.12} has a sense with $h$ replacing $H$. However, we
have conducted the above proof so, because we have in mind to 
generalize our results for a $\Cal C^{1,\alpha}$ mapping,
{\it see} \S7.11 below.
\endremark

\head 
\S7. Relative position of the neighbouring Segre varieties 
\endhead

\subhead 7.1.~Intersection of Segre varieties \endsubhead
We are now in position to complete the proof of Theorem~1.8. Using
Lemma~5.17 and Lemma~6.7, it remains to show that our functions
$\Theta_{w_{q'}',x_{q'}',\beta}'$ extend holomorphically at $0$, for
$\gamma_{w_q,x_q}$ chosen conveniently. For this choice, we are led to
the following dichotomy\,: either $S_{\bar 0}\cap M_1^-=\emptyset$ in
a sufficiently small neighborhood of $0$ or there exists a sequence
$(q_k)_{k\in\N}$ of points of $S_{\bar 0}\cap M_1^-$ tending towards
$0$. In the first case, we shall distinguish two sub-cases. Either
$S_{\bar 0}$ lies under $M_1^-$ or it lies above $M_1^-$. Let us write
this more precisely. We can choose a $\Cal C^\infty$-smooth
hypersurface $H_1$ transverse to $M$ at $0$ with $H_1$ satisfying
$H_1\cap M=M_1$ and $H_1^-\cap M=M_1^-$. Thus $H_1$ together with $M$
divides $\C^n$ near $0$ in four connected parts. We wanted to say that
either $S_{\bar 0}\cap H_1^-$ is contained in the lower left quadrant
$H_1^-\cap M^-$ or it is contained in the upper left quadrant
$H_1^-\cap M^+$. To summarize, we have distinguished three cases\,:
{\it Case~I}\,: $S_{\bar 0}\cap M_1^-\neq \emptyset$ in every
neighborhood of $0$\,; {\it Case~II}\,: $S_{\bar 0}\cap H_1^-\subset
M^-$\,; {\it Case~III}\,: $S_{\bar 0}\cap H_1^-\subset M^+$.  In the
first two cases, the Segre variety
$S_{\overline{\gamma_{w_q,x_q}(0)}}$ will intersect $D\cup \Omega$ and
the neighborhood $\omega(\Sigma_{\gamma_{w_q,x_q}})$ will always
contain $0$\,: the extension will then be an easy direct application
of the Hanges-Treves extension theorem ({\it see} the details
below). The third case could be a priori the most delicate one. But we
can already delineate the following crucial geometric property, which
says that Lemma~5.17 will apply.

\proclaim{Lemma~7.2}
If $S_{\bar 0} \cap H_1^-$ is contained in $M^+$, then $0$ lies in the
lower side $\Sigma_{\gamma_{w_q,x_q}}^-$ for every arc
$\gamma_{w_q,x_q}\subset \subset M_1^-$.
\endproclaim

\demo{Proof}
The real equation of $M$ is given by $y=h(w,\bar w,x)$, where $h$ is a
certain converging real power series satisfying $h(0)=0$, $dh(0)=0$
and $h(w,0,x)\equiv 0$. We can assume that the ``minus'' side $D=M^-$
of automatic extension of CR functions is given by $\{y<h(w,\bar
w,x)\}$. Replacing $x$ by $(z+\bar z)/2$ and $y$ by $(z-\bar z)/2i$,
and solving with respect to $z$, we get for $M$ an equation as above,
say $z=\bar z+i \, \bar\Xi(w,\bar t)$ (we have $\bar \Theta(w,\bar
t)\equiv \bar z+i \, \bar\Xi(w,\bar t)$ in our previous notations),
with $\bar\Xi(0,\bar t)\equiv 0$. Clearly, every such arc
$\gamma_{w_q,x_q}$ contains a point $p\in M$ of coordinates
$(w_p,0+ih(w_p,\bar w_p,0))$ with $u_{1,p}<0$ (indeed, by
construction, these arcs are all elongated and almost directed along
the $x$-coordinate lines since $dh(0)=\text{\rm Id}$). By assumption,
we have $h(w_p,\bar w_p,0)<0$. Then the Segre variety $S_{\bar p}$
(which is a leaf of $\Sigma_{\gamma_{w_q,x_q}}$), has the equation
$z=-ih(w_p,\bar w_p,0)+i \, \bar\Xi(w,\bar w_p,-ih(w_p,\bar w_p,0))$.
Therefore, the intersection point $\{w=0\}\cap S_{\bar p}\subset
\Sigma_{\gamma_{w_q,x_q}}$ whose
coordinates are $(0,-ih(w_p,\bar w_p,0))$ clearly lies over the
origin, which completes the proof of Lemma~7.2.
\qed
\enddemo

\subhead 7.3.~Extension across $(M,0)$ of the components
$\Theta_{w_{q'}',x_{q'}',\beta}'$ \endsubhead We are now prepared to 
complete the proof of Theorem~1.8.
We first choose $\delta, \, \eta, \, \varepsilon$ and various $\v w_q
\v, \v x_q\v<\varepsilon$ as in Lemma~5.17 and we consider the
associated arcs $\gamma_{w_q,x_q}$, $\gamma_{w_{q'}',x_{q'}'}'$, the
associated mapping $h_{w_{q'}',x_{q'}'}$ and the associated reflection
function $\Cal R_{w_{q'}',x_{q'}', h_{w_{q'}',x_{q'}'}}'$. By
Lemmas~6.7 and ~5.17, for each such choice of $(w_q,x_q)$, then all
the components $\Theta_{w_{q'}',x_{q'}',\beta}'$ extend
holomorphically to $D\cup [\Sigma_{\gamma_{w_q,x_q}}^-\cap
\Delta_n(0,\eta)]$. Our goal is to show that for suitably chosen
$\gamma_{w_q,x_q}$ in Cases I, II and III, then the components 
$\Theta_{w_{q'}',x_{q'}',\beta}'$ extend holomorphically across 
$(M,0)$ (afterwards, in Lemma~\S7.10 below, we shall establish 
the desired final Cauchy estimate).

\subhead 7.4.~Case I\endsubhead
In Case~I, we choose an arc $\gamma_{w_q,x_q}$ passing through one of
the points $q_k\in M\cap S_{\bar 0}$ sufficiently close to $0$, with
$\v q_k\v < \eta/2$ and with $\v w_q\v, \v x_q\v <\varepsilon$. Since
$q_k\in M_1^-\cap S_{\bar 0}$, we have $0\in S_{\bar q_k}$. However,
it certainly can happen that $q_k$ belongs to $E_{M'}'$ and then
Lemma~6.7 fails to apply. Fortunately, we can shift slightly
$\gamma_{w_q,x_q}$ in order that $\gamma_{w_{q'}',x_{q'}'}'\cap
E_{M'}'=\emptyset$ (remember\,:
$h(\gamma_{w_q,x_q})=\gamma_{w_{q'}',x_{q'}'}'$). 
Since the Hausdorff dimension of $E_{M'}'$ does
not exceed $2n-3$, it is clear that we can choose such an arc
$\gamma_{w_q,x_q}$ with $q_k$ arbitrarily to close $\gamma_{w_q,x_q}$
and $\gamma_{w_{q'}',x_{q'}'}\cap E_{M'}'$ empty. We let $\tilde{q}$
denote a point of $\gamma_{w_q,x_q}$ which is the nearest to $q_k$. We
have $\v q_k- \tilde{q}\v<< \varepsilon, \, \eta, \, \delta$. By
Lemma~4.5, the components of the reflection function associated with
our fixed coordinate system vanishing at $0$ ({\it i.e.} for
$(w_{q'}',x_{q'}')=(0,0)$) all converge in the polydisc
$\Delta_n(\tilde{q},r(\tilde{q}))$ and they satisfy a Cauchy estimate
there. This implies that the components of the reflection function
associated with $\gamma_{w_{q'}',x_{q'}'}'$ converge in the polydisc
$\Delta_n(q_k,r(q_k)/2)$ and satisfy a Cauchy estimate there
(remember\,: $q\mapsto r(q)$ is lower semi-continuous). To summarize,
we have got\,:

\proclaim{Lemma~7.5}
There exist arcs $\gamma_{w_q,x_q}$ passing through a point
$\tilde{q}$ arbitrarily close to $q_k$ such that the associated
components $\Theta_{w_{q'}',x_{q'}',\beta}'$ extend holomorphically to
$$
D\cup [\Sigma_{\gamma_{w_q,x_q}}^-\cap
\Delta_n(0,\eta)] \cup \Delta_n(q_k,r(q_k)/2).
\tag 7.6
$$
\endproclaim
\noindent
Applying then Lemma~6.9, we see that the neighborhood
$\omega(S_{\bar{\tilde{q}}})$ to which holomorphic extension holds
will contain $0$ if such a $\tilde{q}$ is sufficiently close to $q_k$
(remember\,: $0\in S_{\bar q_k}$). Further, choosing this neighborhood
$\omega(S_{\bar{\tilde{q}}})$ to be of nice tubular form (shrinking it
a bit if necessary), we can certainly assure that its intersection
with the open set~\thetag{7.6} is connected. Case I is done.

\subhead 7.7.~Case II\endsubhead
Case~II, is treated almost the same way. Since $S_{\bar 0}\cap H_1^-$
is contained in $D$, we can choose a fixed point $\tilde{q}$ of
$S_{\bar 0}$ which belongs to $\Delta_n(0,\eta/2)$. Of course, there
exists a positive radius $\tilde{r}>0$ such that the polydisc
$\Delta_n(\tilde{q},\tilde{r})$ is contained in $D$ in which the
components $\Theta_{w_{q'}',x_{q'}',\beta}'$ satisfy a Cauchy estimate, for
all $\v w_{q'}'\v, \v x_{q'}'\v$ sufficiently small. We thus come down to a 
situation similar to that of Lemma~7.5. Case~II is done.

\subhead 7.8.~Case III\endsubhead
For Case~III, thanks to Lemma~7.2, we know already that $0$ belongs to
the lower side $\Sigma_{\gamma_{w_q,x_q}}^-$. Thus Case~III follows
immediately from the application of Lemmas~6.7 and~5.17 summarized in
\S7.3 above.  Case~III is done.

\subhead 7.9.~Extension across $M$ of the reflection function 
\endsubhead To deduce that the
reflection function extends holomorphically at $0$, it remains to
establissh a final Cauchy estimate.

\proclaim{Lemma~7.10}
Let $\v w_{q'}'\v, \v x_{q'}'\v < \varepsilon$ and assume that the
components $\Theta_{w_{q'}',x_{q'}',\beta}'(h(t))$ extend as
holomorphic functions $\varphi_{w_{q'}',x_{q'}',\beta}'(t)$ at
$0$. Then there exist constants $C_0>0$, $r_0>0$ such that $\v
\varphi_{w_{q'}',x_{q'}',\beta}'(t)\v < C_0^{\v \beta\v+1}$ for $\v
t\v < r_0$.
\endproclaim

\demo{Proof}
The formal power series $h_{\Cal F}(t)$ at $0$ satisfies
$\Theta_{w_{q'}',x_{q'}',\beta}'(h_{\Cal F}(t))\equiv
\varphi_{w_{q'}',x_{q'}',\beta}'(t)$. Thanks to the Artin approximation
theorem, there exists a converging power series $H(t)$ such that
$\Theta_{w_{q'}',x_{q'}',\beta}'(H(t))\equiv
\varphi_{w_{q'}',x_{q'}',\beta}'(t)$. Then the easy Cauchy estimate for
the composition of two holomorphic functions yields Lemma~7.10. The 
proof of Theorem~1.8 is complete.
\qed
\enddemo

\subhead 7.11.~Mappings of lesser regularity\endsubhead
Clearly, our constructions in \S2,3,4,5,6,7 are valid for a $\Cal
C^{1,\alpha}$ ($0<\alpha<1$) CR diffeomorphism $h$, except the last
step, namely except the proof of the Cauchy estimates for our
functions $\varphi_{w_{q'}',x_{q'}',\beta}'(t)$, whose holomorphicity
is established along the same lines (for the $\Cal C^{1,\alpha}$
version of the Hanges-Treves theorem, one has to use the work of
Tumanov [Tu2]). However, in the holomorphically nondegenerate case,
we can complete the following\,:

\proclaim{Theorem~7.12}
If $h$ is a $\Cal C^{1,\alpha}$-smooth CR-diffeomorphism, if $(M,p)$
is minimal and if $(M',p)$ is holomorphically nondegenerate, then the
reflection function $\Cal R_h'$ extends holomorphically to a
neighborhood of $p\times \overline{p'}$.
\endproclaim

\noindent
(However, we are still unable to conclude that $h$ is real analytic.)
Indeed, first of all, it is not difficult to deduce then that such a $\Cal
C^{1,\alpha}$ map extends as a correspondence at the point $p_1\in
M_1$ (by taking the uniquely defined irreducible component of $\Cal
C_h'$ which contains the graph of $h$). Consequently, we get
polynomial identities with coefficients holomorphic in a neighborhood
of $p_1$ for the components of the mapping $h$. Finally, to prove that
the reflection function extends holomorphically at $p_1$, namely to
get the desired Cauchy estimates, we can then apply the same scheme of
proof as for instance in [BJT] or [DP].

\vfill
\enddocument
\end